\documentclass[preprint,12pt]{elsarticle}

\usepackage[utf8]{inputenc}
\usepackage[a4paper,left=3cm,right=3cm,top=2.5cm,bottom=2.5cm]{geometry}

\usepackage{tensor, bm}
\usepackage{amssymb,amsmath}

\usepackage{algorithm}
\usepackage{algpseudocode}

\usepackage{lineno}
\usepackage{pgfplots}
\pgfplotsset{compat=1.17}
\usepackage{pstricks,pst-node, pst-fun}
\usepackage{dsfont}
\usepackage[spanish,USenglish,american]{babel}

\usepackage{siunitx}
\usepackage{arydshln}
\usepackage{relsize}
\usepackage{setspace}
\usepackage{caption,subcaption}
\usepackage[hidelinks]{hyperref}
\usepackage[symbol]{footmisc}

\usepackage{tikz}
\usetikzlibrary{shapes.geometric, arrows}
\tikzstyle{startstop} = [rectangle, rounded corners, minimum width=3cm, minimum height=1cm,text centered, draw=black, fill=red!30]
\tikzstyle{io} = [trapezium, trapezium left angle=70, trapezium right angle=110, minimum width=3cm, minimum height=1cm, text centered, draw=black, fill=yellow!30] 
\tikzstyle{process} = [rectangle, minimum width=3cm, minimum height=1cm, text centered, draw=black, fill=blue!30]
\tikzstyle{decision} = [diamond, minimum width=1cm, minimum height=1cm, text centered, draw=black, fill=green!30]

\graphicspath{{./figures/}}

\usepackage{easybmat}
\usepackage{array} 
\usepackage{booktabs,colortbl}	     

\definecolor{Gray}{gray}{0.9}

\newcolumntype{C}[1]{>{\centering\arraybackslash}m{#1}}
\newcolumntype{I}{!{\vrule width 1pt}}

\newcommand{\W}{\mathcal{W}}
\newcommand{\h}{\mathcal{H}}
\newcommand{\heps}{\mathcal{H}^{\bm \varepsilon}}
\newcommand{\hsig}{\mathcal{H}^{\bm \sigma}}

\makeatletter 
\def\ps@pprintTitle{%
  \let\@oddhead\@empty
  \let\@evenhead\@empty
  \def\@oddfoot{\centerline{\thepage}}%
  \let\@evenfoot\@empty
}
\makeatother

\begin{document}

\begin{frontmatter}



\title{An elastic properties-based topology optimization algorithm for linear orthotropic, functionally graded materials}


\renewcommand{\thefootnote}{\fnsymbol{footnote}}

\author{Ismael Ben-Yelun$^{a,}$\footnote[1]{Corresponding author\\ \hspace*{6mm}%
\textit{Email address}: \texttt{i.binsenser@upm.es} (I. Ben-Yelun).}, Víctor Riera$^a$, Luis Saucedo-Mora$^{a,b,c}$, Miguel Ángel Sanz$^a$, Francisco Javier Montáns$^{a,d}$}

\address{$^a$ E.T.S. de Ingeniería Aeronáutica y del Espacio, Universidad Politécnica de Madrid, Pza. Cardenal Cisneros 3, 28040, Madrid, Spain\\

$^b$ Department of Materials, University of Oxford, Parks Road, Oxford, OX1 3PJ, UK\\

$^c$ Department of Nuclear Science and Engineering, Massachusetts Institute of Technology,  MA02139, USA\\

$^d$ Department of Mechanical and Aerospace Engineering, Herbert Wertheim College of Engineering, University of Florida, FL32611, USA}

\begin{abstract}
Topology optimization (TO) has experienced a dramatic development over the last decades aided by the arising of metamaterials and additive manufacturing (AM) techniques, and it is intended to achieve the current and future challenges. In this paper we propose an extension for linear orthotropic materials of a three-dimensional TO algorithm which directly operates on the six elastic properties -- three longitudinal and shear moduli, having fixed three Poisson ratios -- of the finite element (FE) discretization of certain analysis domain. By performing a gradient-descent-alike optimization on these properties, the standard deviation of a strain-energy measurement is minimized, thus coming up with optimized, strain-homogenized structures with variable longitudinal and shear stiffness in their different material directions. To this end, an orthotropic formulation with two approaches -- direct or strain-based and complementary or stress-based -- has been developed for this optimization problem, being the stress-based more efficient as previous works on this topic have shown.

The key advantages that we propose are: (1) the use of orthotropic ahead of isotropic materials, which enables a more versatile optimization process since the design space is increased by six times, and (2) no constraint needs to be imposed (such as maximum volume) in contrast to other methods widely used in this field such as Solid Isotropic Material with Penalization (SIMP), all of this by setting one unique hyper-parameter. Results of four designed load cases show that this orthotropic-TO algorithm outperforms the isotropic case, both for the similar algorithm from which this is an extension and for a SIMP run in a FE commercial software, presenting a comparable computational cost. We remark that it works particularly effectively on pure shear or shear-governed problems such as torsion loading.
\end{abstract}

\begin{keyword}
Topology optimization \sep Orthotropic materials \sep Functionally graded \sep Mechanical metamaterials
\end{keyword}

\end{frontmatter}


\section{Introduction}

Novel and more flexible methods for topology optimization (TO) are becoming increasingly popular, leveraged by the advances in additive manufacturing \cite{brackett2011topology} and the arise of metamaterials \cite{diaz2010topology, bertoldi2017flexible}---the latter preceded by the greater computational capacity of modern computers. TO is an area within the structural optimization wherein a given domain $\Omega$ subjected to certain boundary conditions $\Gamma_N \cup \Gamma_D$ is sought to present an optimized final shape by minimizing one of its macroscopic features e.g. mean compliance or total mass. It has been historically addressed with a Finite Element (FE) approach, therefore similar procedures have been followed herein \cite{bathe2006finite}.

In the search for more versatile techniques, the use of orthotropic materials seems a promising choice due to their characteristics \cite{jones2018mechanics} e.g. high stiffness and strength-to-weight ratio \cite{karatacs2018review}. An example of their usage might be seen in the aerospace sector, where more than 50\% of the primary structure of the A350 (latest Airbus aircraft programme) is made of Carbon Fiber Reinforced Polymer (CFRP) laminates i.e. aligned carbon fibers embedded in epoxy matrix \cite{mcilhagger2020manufacturing}. In fact, this links it with the TO in the sense that originally, the practical applications of this subject were mainly carried out in this industry \cite{suzuki1991homogenization}---a review of several aeronautical utilizations has been recently done by Zhu et al. \cite{zhu2016topology}. Regarding the mechanical properties, constructing different material architectures enable a wider range of possibilities. Through the design of orthotropic materials (e.g. composites), one may achieve a purposeful anisotropy in the elastic properties \cite{pedersen1987sensitivity, pedersen1989optimal}. A more recent example is the development of metamaterials, wherein macro-scale behaviour is controlled by tuning the design parameters at micro, unit-cell level \cite{bertoldi2017flexible}, allowing the modification of elastic properties e.g. Poisson ratio within a certain physical range.

This paper belongs to the area of structural optimization, which several authors have subdivided into three: size optimization, shape optimization and topology optimization \cite{christensen2008introduction, tsavdaridis2019application}. While size optimization is clearly identified and separated (e.g. obtaining the optimal thicknesses of a set of parts in a structural component that satisfy given constraints), the other two are devoted to respectively optimize the boundary $\partial \Omega$ (i.e. shape), and the connectivity within the analysis domain $\Omega$ (i.e. topology). Since these two last may overlap to some extent, another classification by methods has been made by Deaton and Grandhi \cite{deaton2014survey, yago2021topology}, comprising of three groups. The first one is homogenization (e.g. density-based), which began with the pioneering works of Bends{\o}e resulting in the commonly-used Solid Isotropic Material with Penalization (SIMP) method \cite{Bendse1989, Bendsoe1999}. This is the method against which the proposed algorithm is compared in this paper from a conservative point of view. In the second group  lie the evolutionary methods (or hard-kill), such as the Evolutionary Structural Optimization (ESO) proposed by Xie and Steven \cite{xie1993simple} or its birectional extension (BESO), by Huang and Xie \cite{querin1998evolutionary}. More specifically, Xie and co-workers developed and algorithm wherein the FE removal is carried out based on their contribution to the total compliance \cite{huang2008topology, huang2009bi}. In that sense, this work implements a similar feature since elements are removed from the domain attending to their energy contribution and their derivatives (with respect to elastic properties) to the whole domain. Finally, level-set (boundary variation) methods attempt to fix the artifacts that may arise in homogenization methods by defining level-set curves or surfaces that deal better with more complex topologies \cite{wu2017level}. These previous methods were initially meant for isotropic, linear materials, but other effects such as geometrical or material non-linearities have been included in TO algorithms as well. Regarding material non-linearities, Yuge and Kikuchi addressed an elasto-plastic analysis \cite{yuge1995optimization}, Bends{\o}e et al. studied softening material \cite{bendsoe1996optimization} and Zhang et al. addressed multi-material TO considering material non-linearities \cite{zhang2018multi}, among others. Geometric non-linearities were also considered by Pedersen and Sigmund in their large-displacement TO work \cite{pedersen2001topology}, as well as the numerical methods developed by Bruns et al. to address a non-linear elasticity problem such as snap-through \cite{bruns2002numerical}.

A linear extension of material behavior that may be posed consists of considering linear orthotropic materials---different orthotropic applications for TO might be found in the literature. However, these studies have nearly always been limited to bi-dimensional, plane problems comprised of unit cells and their orientation in space. These unit cells are virtually a metamaterial configuration: isotropic material elements with a void at their centre, thus presenting a global anisotropic behavior, in overall. Therefore, those are optimal orientation problems, and the first published works on the subject were carried out by Pedersen \cite{pedersen1989optimal, pedersen1990bounds, pedersen1991thickness}, who developed the models within a strains-based framework. Suzuki and Kikuchi \cite{suzuki1991homogenization}, and Díaz and Bendsoe \cite{diaz1992shape} -- whose works were devoted to the optimal orientation of the so-called unit cells in plane-stress problems applied for TO algorithms as well -- highlighted that the stress-based formulation is more efficient. Cheng, Kikuchi and Ma addressed similarly this topic from the point of view of Optimal Material Distribution (OMD) \cite{cheng1994improved}, also agreeing with the efficiency of the stress-based formulation for these orthotropic TO problems. 

Regarding the development of TO methods for orthotropic materials in the recent years, research has been likewise conducted yet it is still scarce. Gea and Luo provided closed-form solutions for both strain-based and stress-based orthotropic TO problems modelled through these plane unit-cells \cite{gea2004stress}. Jia et al. continued with this approach \cite{jia2008topology}, presenting a SIMP-alike model optimizing at the same time the density of FEs and orientation of the cells. On the other side, Luo and Gea also proposed a model based on plate structures, similarly to laminates for composites materials \cite{luo1998optimal}, and Stegmann and Lund developed an optimization algorithm for laminates, thus belonging to Discrete Material Optimization (DMO) methods \cite{stegmann2005discrete}. We remark that our study is not limited to only a discrete amount of plies (e.g. laminate) but to all kind of orthotropic materials, thus enabling the availability of continuous derivatives -- which are useful for gradient optimization methods -- or the usage of this method in conjunction with AM techniques or metamaterials. More recently, Page et al. \cite{page2016topology} have proposed a similar method insofar a volume restriction is imposed, but applied to heat transfer problems. Another application with AM techniques is proposed by Li et al. \cite{lee2021design}, also using these unit cells with orthotropic behaviour. A SIMP scheme is used likewise, with the previously mentioned drawback of having to impose a volume restriction.

In this paper, we propose a novel TO algorithm suitable for linear orthotropic materials to come up with three-dimensional optimized structures with minimum local compliance through two different approaches: strain-based, i.e. using the strain energy function, and stress-based, i.e. using the complementary strain (or \textit{stress}) energy function. This optimization algorithm updates at once the elastic properties (the 6 stiffness moduli, longitudinal and shear, letting the three Poisson ratios fixed) of each FE in which the domain is discretized. These updates are addressed by linearizing a decoupled form of the energy function, using their derivatives with respect to said elastic properties to perform gradient-descent-alike update step. Thus, by operating directly on the elastic properties, the step of relating them to a intermediate variable e.g. density or mass is saved. Our model presents particularities typical of orthotropic materials, for instance, the coupling of longitudinal moduli which appears in the terms of the volumetric part of the constitutive matrix in the strain-based approach. Hence the need of introducing the complementary formulation, in an uncoupled fashion preferably. For this purpose, a similar decoupling of the (complementary) strain energy function as the one performed in Amores et al. is followed \cite{amores2021finite}. This formulation is proved to be more efficient, as former research have demonstrated likewise \cite{suzuki1991homogenization, diaz1992shape}.

Therefore, there are two key features that we have developed: The first one is the implementation of a numerical formulation for orthotropic materials in TO framework, leading to a more versatile optimization due to the extension of the design space---by six times with respect to the isotropic material case. The second one is that a volume constraint is no longer needed as in other classical methods. Similarly, heuristic sensitivity filters are not required. Additionally, since the final goal is to achieve structures with functionally graded properties, no penalization is applied. All of this requires fixing only one step-update hyperparameter needed for the optimization. In order to assess the performance of both methods, four examples are studied. These are designed such that the nature of different load cases (tensile-compressive, pure shear, combination of both) are represented. Then, the results of the algorithm are compared with their isotropic analogue \cite{saucedo2023updated, ben2023topology} and with simulation runs carried out in the commercial FEM-CAE software \textsc{OptiStruct}, from \textsc{Altair} \cite{optistruct}.

This paper is organized as follows. First, in Section~\ref{Sec:theoretical_framework} the theoretical framework is established, in which the formulation proposed is developed. Then, in Section~\ref{Sec:methodology} the methodology and particularities of this method are outlined, considering important restrictions and sketching the most important algorithms. Finally, in Section~\ref{Sec:results} the previously commented results are displayed and analysed, concluding with a number of final remarks in Section~\ref{Sec:concluding_remarks}.

\section{Theoretical framework}
\label{Sec:theoretical_framework}

In contrast to other TO methods, the objective function to minimize is the standard deviation $s_{\h^{\bm \varepsilon}}$ of a strain-energy like variable $\h^{\bm \varepsilon}$, where $\bm \varepsilon$ stands for strains---see Saucedo et al. \cite{saucedo2023updated} and Ben-Yelun et al. \cite{ben2023topology} for more details. This is achieved by directly operating on the different stiffness of the material, which enables the elimination of intermediate variables (e.g. density, as it is done in density-based topology optimization methods such as SIMP \cite{Bendse1989, Bendsoe1999}).

The novelty with respect to other works is the nature of the considered material: orthotropic. A linear orthotropic material contains 9 elastic properties: Young moduli on 3 principal material directions $E_1$, $E_2$, $E_3$, 3 shear moduli referred to the same directions $G_{12}$, $G_{13}$, $G_{23}$ and three of the six Poisson coefficients e.g. $\nu_{12},\nu_{13},\nu_{23}$---the remaining three Poisson coefficients are determined due to tensor symmetry i.e. $\nu_{ij}/E_i = \nu_{ji}/E_j$. Thus, in terms of the optimization problems, this implies a wider range in design space compared to the linear isotropic case, which only presents two elastic parameters, for instance. Therefore, since the exploration of more possibilities is allowed, optimal designs might outperform the baseline results provided by its isotropic counterpart.

\subsection{Energy split into volumetric and deviatoric parts}

First of all, the volumetric-deviatoric split on the (total) energy is performed so that

\begin{equation}
\mathcal{W}^{\square} := \mathcal{W}^{\square,v} + \mathcal{W}^{\square,d},
\end{equation}
where, for the sake of compactness, $\square = \left\lbrace \bm \varepsilon, \bm\sigma\right\rbrace$ can represent either direct (strains-based) or complementary (stress-based) energy, respectively. With that, the first term $\mathcal{W}^{\square, v}$ contains information about the longitudinal stiffness and their coupling (i.e. Young moduli and Poisson ratio) whilst the shear stiffness is involved in the second term $\mathcal{W}^{\square, d}$. It is important to remark that if an isotropic material optimization with fixed Poisson ratio is chosen, this splitting would be futile since there is only one variable to optimize: the Young modulus.

The field variables of this structural problem are obtained through finite-element modeling (FEM). Thus, the domain is discretized into $n_e$ elements, hence the energy might be expressed as the sum of the contributions of each element:

\begin{equation}
\mathcal{W}^{\square} = \sum_{e=1}^{n_e}\mathcal{W}_e^{\square} = \sum_{e=1}^{n_e} \W_e^{\square, v} + \sum_{e=1}^{n_e} \W_e^{\square,d}.
\end{equation}

Assuming that the material works within its linear elastic regime, both direct and complementary elastic energies might be expanded through stiffness and compliance matrices---$\mathbf{D}$ and $\mathbf{D}^{-1}$ in Voigt notation, respectively. For materials with longitudinal-shear decoupling (as in the case of isotropic and orthotropic materials, among others) these matrices could be separated as follow into the following submatrices, being null the part relating extension and shear

\begin{equation}
\left[\mathbf{D}\right] =
\left[
\begin{BMAT}{c;c}{c;c}
\left[\mathbf{D}^v\right] & \left[\bm{0}\right] \\ 
\left[\bm{0}\right] & \left[\mathbf{D}^d\right]
\end{BMAT}
\right], 
\quad
\left[\mathbf{D}^{-1}\right] =
\left[
\begin{BMAT}{c;c}{c;c}
\left[\left(\mathbf{D}^v\right)^{-1}\right] & \left[\bm{0}\right] \\
\left[\bm{0}\right] & \left[(\mathbf{D}^d )^{-1}\right]
\end{BMAT}
\right],
\end{equation}
with \cite{chaves2014mecanica}

\begin{equation}
\begin{array}{c}
[\mathbf{D}^v]= \dfrac{1}{\chi}
\begin{bmatrix}
E_{1} \left(1-\nu_{23}^{2} \frac{E_{3}}{E_{2}}\right) &
\nu_{12} E_{2}+\nu_{23} \nu_{13} E_{3} &
\nu_{13} E_{3}+\nu_{12} \nu_{23} E_{3}\\
& E_{2} \left(1-\nu_{13}^{2} \frac{E_{3}}{E_{1}}\right) &
\nu_{23} E_{3}+\nu_{12} \nu_{13} \frac{E_{2} E_{3}}{E_{1}}\\
\rm{(sym)} & & E_{3} \left(1-\nu_{12}^{2} \frac{E_{2}}{E_{1}}\right)\\
\end{bmatrix},
\\[7.5ex]
\chi = 1 - \dfrac{\nu_{12}^2E_2}{E_1} - \dfrac{\nu_{23}^2E_3}{E_2} - 2\nu_{12}\nu_{23}\nu_{13}\dfrac{E_3}{E_1},
\end{array}
\label{Eq:2_D_v}
\end{equation}

\begin{equation}
[(\mathbf{D}^v)^{-1}] =
\begin{bmatrix}
1/E_1 & -\nu_{21}/E_2 & -\nu_{31}/E_3 \\
-\nu_{12} / E_1 & 1/E_2 & -\nu_{32}/E_3 \\
 -\nu_{13} / E_1 & -\nu_{23} / E_2 & 1 / E_3
\end{bmatrix},
\label{Eq:2_D_d}
\end{equation}
and
\begin{equation}
[\mathbf{D}^d] = \begin{bmatrix}
G_{12} & 0 & 0 \\ 0 & G_{13} & 0 \\ 0 & 0 & G_{23}
\end{bmatrix},
\quad
[(\mathbf{D}^d)^{-1}] = \begin{bmatrix}
1/G_{12} & 0 & 0 \\ 0 & 1/G_{13} & 0 \\ 0 & 0 & 1/G_{23}
\end{bmatrix}.
\label{Eq:2_inverse_D}
\end{equation}

In the following subsections, the optimization algorithms developed for the direct and complementary formulation are addressed.

\subsection{Direct (strain-based) approach with separated influence of elastic properties}
\label{Subsec:2-direct}

In this first approach, the elastic energy of the $e-$th element $\W_e^{\bm \varepsilon}$ might be expanded as follows by performing numerical integration

\begin{equation}
\W_e^{\bm \varepsilon} =  \dfrac{1}{2}\sum_{p=1}^{n_p}\underline{\bm{\varepsilon}}\left(\mathbf{x}_p\right)\cdot \mathbf{D}_{k}^v \underline{\bm{\varepsilon}}\left(\mathbf{x}_p\right) J_p w_p
+ 
\dfrac{1}{2}\sum_{p=1}^{n_p}\underline{\bm{\gamma}}\left(\mathbf{x}_p\right)\cdot \mathbf{D}_{k}^d \underline{\bm{\gamma}}\left(\mathbf{x}_p\right) J_p w_p,
\label{Eq:2_W_e}
\end{equation}
where $\underline{\bm{\varepsilon}}$ is a vector containing the longitudinal strains along the material directions i.e. $\underline{\bm{\varepsilon}} = [\varepsilon_{1}, \varepsilon_{2}, \varepsilon_{3}]^T$, and $\underline{\bm{\gamma}}$ is the vector which contains the shear strains of the material directions $\underline{\bm{\gamma}} = [\gamma_{12}, \gamma_{13}, \gamma_{23}]^T$. Both $\underline{\bm{\varepsilon}}$ and $\underline{\bm{\gamma}}$ are evaluated in the integration point $p$, as well as the Jacobian of the element transformation to normalized coordinates $J_p$, being their corresponding quadrature weights $w_p$. This is a generalization for elements with $n_p$ integration points.

A similar procedure to the isotropic case might be followed. In these materials, the energy contribution of an element $\W_e^{\bm \varepsilon}$ might be expressed as a product of its Young modulus $E_e$ and an \textit{ad-hoc} variable $\h_e^{\bm \varepsilon}$ which contains all the energy information except the Young modulus i.e. $\W_e^{\bm \varepsilon} = E_e \h_e^{\bm \varepsilon} \Omega_e$ (being $\Omega_e$ the volume of the element $e$), in such a way that $\partial \h_e^{\bm \varepsilon} / \partial E_e = 0$, considering the rest of the variables fixed e.g. strain/stress distributions. The introduction of the volume $\Omega_e$ in the last expression is to prevent the algorithm from being sensitive to different element sizes. A volume-averaged strain energy density variable $\bar{\Psi}^{\bm \varepsilon}$ represents a more suitable alternative, which is defined as

\begin{equation}
\bar{\Psi}_e^{\bm \varepsilon} := \dfrac{1}{\Omega_e}\int_{\Omega_e}\Psi^{\bm \varepsilon}\left(\mathbf{x}\right)\mathrm{d}\Omega_e = \dfrac{\W_e^{\bm \varepsilon}}{\Omega_e}.
\end{equation}

In this way, the stiffening of larger elements by the mere fact of having larger volume than others is avoided. Therefore, it is desirable to model orthotropic materials in a similar way. To this end, a split of the elastic energy on six terms is carried out, each of them containing one of the six moduli in linear orthotropic materials. Ideally, these terms depend explicitly and solely on that elastic property that represent, so that a decoupled form of the energy might be achieved, where in each term a separation of variables with the elastic modulus can be identified. Taking this into account, the energy of each element $e$ might be expressed as

\begin{equation}
\mathcal{W}_e^{\bm \varepsilon} = 
\sum_{i=1}^3\mathcal{W}_{i,e}^{\bm \varepsilon, v} + \sum_{i=1}^3\sum_{j>i}^3\mathcal{W}_{ij,e}^{\bm \varepsilon, d} =
\sum_{i=1}^3\mathcal{W}^{\bm \varepsilon, d}_{i,e}(E_{i,e}) + \sum_{i=1}^3\sum_{j>i}^3\mathcal{W}^{\bm \varepsilon, d}_{ij,e} (G_{ij,e}),
\label{Eq:2_energy_of_an_element}
\end{equation}
where $\mathcal{W}_{i,e}^v$ represents the elastic energy contribution which involves stress and strain in the principal direction $i$ (i.e. the longitudinal stiffness and the Poisson effects), and $\mathcal{W}^d_{ij,e}$ is the contribution of the shear part which operates in directions $i$ and $j$.

Thus, splitting Equation~(\ref{Eq:2_W_e}) in three longitudinal and another three shear parts,

\begin{equation}
\W_e^{\bm \varepsilon} = \dfrac{1}{2}\sum_{i=1}^3\sum_{p=1}^{n_p} \underline{\bm{\varepsilon}}\left(\mathbf{x}_p\right)\cdot \mathbf{D}_{i,e}^v \underline{\bm{\varepsilon}}\left(\mathbf{x}_p\right) J_p w_p
+ 
\dfrac{1}{2}\sum_{i=1}^3 \sum_{j>i}^3 \sum_{p=1}^{n_p}\underline{\bm{\gamma}}\left(\mathbf{x}_p\right)\cdot \mathbf{D}_{ij,e}^d \underline{\bm{\gamma}}\left(\mathbf{x}_p\right) J_p w_p.
\end{equation}

The disadvantage of this formulation lies in the difficulty of separating $\mathbf{D}_{1,k}^v$, $\mathbf{D}_{2,k}^v$ and $\mathbf{D}_{3,k}^v$ since an expression of $E_1$, $E_2$ and $E_3$ appears non-explicitly in each term of $\mathbf{D}_{k}^v$---see Equation~(\ref{Eq:2_D_v}). Conversely, $\mathbf{D}_{ij,e}^d$ are straightforward to relate with their shear modulus $G_{ij,e}$---see Equation~(\ref{Eq:2_inverse_D}).

We can now introduce the volume-averaged strain energy densities $\bar{\Psi}_e^{\bm \varepsilon, v}$ and $\bar{\Psi}_e^{\bm \varepsilon, d}$ by analogy,

\begin{equation}
\W_e^{\bm \varepsilon} = \sum_{i=1}^3 \bar{\Psi}_{i,e}^{\bm \varepsilon, v} \Omega_e
+ 
\sum_{i=1}^3 \sum_{j>i}^3 \bar{\Psi}_{ij,e}^{\bm \varepsilon, d} \Omega_e,
\end{equation}
where the volume by means of numerical integration is $\Omega_e = \sum_p J_p w_p$. Proceeding in an analogous way to the isotropic case, by defining the strain-variables $\h_{i,e}^{\bm \varepsilon}$ and $\h_{ij,e}^{\bm \varepsilon}$ associated with the volumetric and deviatoric contributions, we propose to rewrite the previous terms in the following way

\begin{equation}
\bar{\Psi}_{i,e}^{\bm \varepsilon, v} = \dfrac{\sum_p \underline{\bm{\varepsilon}}\left(\mathbf{x}_p\right)\cdot \mathbf{D}_{i,e}^v\underline{\bm{\varepsilon}}\left(\mathbf{x}_p\right)J_p w_p}{2\sum_p J_p w_p} =: E_{i,e}\h_{i,e}^{\bm \varepsilon},
\label{eq:separation_variables_direct}
\end{equation}
and

\begin{equation}
\bar{\Psi}_{ij,e}^{\bm \varepsilon, d} = \dfrac{\sum_p \underline{\bm{\gamma}}\left(\mathbf{x}_p\right)\cdot \mathbf{D}_{ij,e}^d\underline{\bm{\gamma}}\left(\mathbf{x}_p\right)J_p w_p}{2\sum_p J_p w_p} =: G_{ij,e}\h_{ij,e}^{\bm \varepsilon, d}.
\end{equation}

Note that we \textit{assume} the explicit separation of variables in the volumetric part of the direct case in Equation~\eqref{eq:separation_variables_direct} into the Young's moduli $E_{i,e}$ and the strain variable $\h^{\bm \varepsilon}_{i,e}$ in order to follow a similar procedure to \cite{ben2023topology}. This separation is however correct in the deviatoric contribution of the energy since $G_{ij,e}$ is a explicit, separable term in $\bar{\Psi}_{ij,e}^{\bm \varepsilon, d}$ i.e., no assumption is required. 

Therefore, the \textit{ad-hoc} variables $\h_{i,e}^{\bm \varepsilon}$ and $\h_{ij,e}^{\bm \varepsilon}$ are defined as

\begin{equation}
\h_{i,e}^{\bm \varepsilon} := \dfrac{\sum_p\underline{\bm{\varepsilon}}\left(\mathbf{x}_p\right)\cdot \hat{\mathbf{D}}_{i,e}^v \underline{\bm{\varepsilon}}\left(\mathbf{x}_p\right) J_p w_p}{2\sum_p J_p w_p},
\end{equation}
and

\begin{equation}
\h_{ij,e}^{\bm \varepsilon} := \dfrac{\sum_p\underline{\bm{\gamma}}\left(\mathbf{x}_p\right)\cdot \hat{\mathbf{D}}_{ij,e}^d \underline{\bm{\gamma}}\left(\mathbf{x}_p\right) J_p w_p}{2\sum_p J_p w_p}.
\end{equation}

Regarding these introduced variables, the following dimensionless matrices are proposed for the volumetric part (in order to compute $\h_{i,e}^{\bm \varepsilon}$):

\begin{equation}
\hat{\mathbf{D}}_{i,e} = \dfrac{\mathbf{D}_e}{3E_i},
\end{equation}
so the volumetric contribution of the elastic energy is recovered when these reduced matrices are sum over the three material directions. As highlighted before, $\hat{\mathbf{D}}_{i,e}^v$ neither depends solely on $E_i$, nor is this dependence explicit i.e. $\partial \h_{i,e}^{\bm \varepsilon} / \partial E_{i,e} \neq 0$. Hence the advantage of a complementary formulation in which the Young moduli $E_i$ appear explicitly in the terms of the (compliance) matrix. This complementary formulation will be covered later on.

However, this explicit separation is indeed achieved in the deviatoric contribution since the matrices $\hat{\mathbf{D}}^d_{ij,e}$ might be easily obtained as

\begin{equation}
\hat{\mathbf{D}}^d_{ij,e} = \dfrac{\mathbf{D}^d_{ij,e}}{G_{ij,e}},
\end{equation}
leading to the following reduced matrices

\begin{equation}
[\hat{\mathbf{D}}^d_{12,k}] = \begin{bmatrix}
1 & 0 & 0 \\ 0 & 0 & 0 \\ 0 & 0 & 0
\end{bmatrix}, \;
[\hat{\mathbf{D}}^d_{13,k}] = \begin{bmatrix}
0 & 0 & 0 \\ 0 & 1 & 0 \\ 0 & 0 & 0
\end{bmatrix}, \;
[\hat{\mathbf{D}}^d_{23,k}] = \begin{bmatrix}
0 & 0 & 0 \\ 0 & 0 & 0 \\ 0 & 0 & 1
\end{bmatrix}.
\label{Eq:2_D_ij_e_reduced}
\end{equation}

Note that these reduced matrices are constant regardless of the element $e$ that they represent, which leads to $\partial \bar{\Psi}^{\bm \varepsilon,d}_{ij,e} / \partial G_{ij,e} = \h_{ij,e}^{\bm \varepsilon}$.

Having computed $\bm\heps_i$ and $\bm \heps_{ij}$, the total energy might be computed by

\begin{equation}
\W_e^{\bm \varepsilon} = \displaystyle\sum_{i=1}^3 E_{i,e}\h_{i,e}^{\bm \varepsilon}\Omega_e
+ 
\displaystyle\sum_{i=1}^3 \sum_{j>i}^3 G_{ij,e}\h_{ij,e}^{\bm \varepsilon} \Omega_e.
\label{Eq:2_W_e_decoupled_and_H}
\end{equation}

\subsection{Complementary (stress-based) approach with separated influence of elastic properties}
\label{Subsec:2-complementary}

In the linear elastic regime, the complementary elastic energy is equal to the direct elastic energy i.e. $\W^{\bm \sigma} = \W^{\bm \varepsilon}$, so the computation is analogous. The complementary strain energy function might be expressed as \cite{holzapfel2002nonlinear}

\begin{equation}
\mathcal{W}^{\bm \sigma} = 
\int_\Omega \dfrac{1}{2}\bm{\sigma}:\bm{\varepsilon}\left(\bm{\sigma}\right)\mathrm{d}\Omega =
\int_\Omega \dfrac{1}{2}\bm{\sigma}:\mathbb{S}:\bm{\sigma}\,\mathrm{d}\Omega
=
\sum_{e=1}^{n_e}\int_{\Omega_e}\dfrac{1}{2}\bm{\sigma}:\mathbb{S}:\bm{\sigma}\,\mathrm{d}\Omega_e =: \sum_{e=1}^{n_e}\W_e^{\bm \sigma},
\label{Eq:2_energy_complementary_definition}
\end{equation}
where $\mathbb{S}$ is the fourth-order compliance tensor, defined such that $\bm{\varepsilon} = \mathbb{S}:\bm{\sigma}$. Thus, the complementary strain energy contribution of element $e$, $\W^{\bm \sigma}_e$, could be defined as

\begin{equation}
\W_e^{\bm \sigma} = \dfrac{1}{2}\sum_{p=1}^{n_p}\underline{\bm{\sigma}}\left(\mathbf{x}_p\right)\cdot (\mathbf{D}_e^v)^{-1} \underline{\bm{\sigma}}\left(\mathbf{x}_p\right) J_p w_p
+ 
\dfrac{1}{2}\sum_{p=1}^{n_p}\underline{\bm{\tau}}\left(\mathbf{x}_p\right)\cdot (\mathbf{D}_e^d)^{-1} \underline{\bm{\tau}}\left(\mathbf{x}_p\right) J_p w_p,
\label{Eq:2_W_e_c}
\end{equation}
where numerical integration with $n_p$ integration points is likewise performed, and vectors $\underline{\bm{\sigma}} = [\sigma_{1}, \sigma_{2}, \sigma_{3}]^T$ and $\underline{\bm{\tau}} = [\tau_{12}, \tau_{13}, \tau_{23}]^T$ are defined analogously to $\underline{\bm{\varepsilon}}$ and $\underline{\bm{\gamma}}$. Note that $\underline{\bm{\sigma}} = (\mathbf{D}^v)^{-1}\underline{\bm{\varepsilon}}$ and $\underline{\bm{\tau}} = (\mathbf{D}^d)^{-1}\underline{\bm{\gamma}}$. Similarly to the direct approach, by splitting the complementary energy function $\mathcal{W}^{\bm \sigma}$ into volumetric and deviatoric contributions and being considered the previously referred finite-element discretization,

\begin{equation}
\W_e^{\bm \sigma} = \dfrac{1}{2}\sum_{i=1}^3 \sum_{p=1}^{n_p}\underline{\bm{\sigma}}\left(\mathbf{x}_p\right)\cdot (\mathbf{D}_{i,e}^v)^{-1} \underline{\bm{\sigma}}\left(\mathbf{x}_p\right) J_p w_p
+ 
\dfrac{1}{2}\sum_{i=1}^3 \sum_{j>i}^3\sum_{p=1}^{n_p} \underline{\bm{\tau}}\left(\mathbf{x}_p\right)\cdot (\mathbf{D}_{ij,e}^d)^{-1} \underline{\bm{\tau}}\left(\mathbf{x}_p\right) J_p w_p.
\label{Eq:2_W_e_c_decoupled}
\end{equation}

In this case, the dependence of the terms of compliance matrices is explicit on the elastic parameters that they are associated with, so a direct identification of $(\mathbf{D}_{i,e}^v)^{-1}$ and $(\mathbf{D}_{ij,e}^d)^{-1}$ from the total matrices might be carried out. This represents an advantage over the direct case which results from a better or less complex formulation of the volumetric part.

Following the analogy with respect to the isotropic case, the complementary formulation in these materials states that the volume-averaged energy density $\bar{\Psi}^{\bm \sigma}_e$ might be expressed as $\bar{\Psi}_e^{\bm \sigma} = \h_e^{\bm \sigma} / E_e$, being defined an \textit{ad-hoc} stress variable $\h_e^{\bm \sigma}$. Thus, it is satisfied that $\partial \h^{\bm \sigma}_e / \partial E_e = 0$ assuming that the stress distribution remains fixed.

Extrapolating it to the orthotropic case, Equation~(\ref{Eq:2_W_e_c_decoupled}) might be rewritten as

\begin{equation}
\W_e^{\bm \sigma} = \sum_{i=1}^3 \bar{\Psi}_{i,e}^{\bm \sigma,v}\Omega_e +
\sum_{i=1}^3\sum_{j>i}^3\bar{\Psi}^{\bm \sigma,d}_{ij,e} \Omega_e,
\end{equation}
which allows us to group terms into the following variables

\begin{equation}
\bar{\Psi}_{i,e}^{\bm \sigma,v}
=
\dfrac{\sum_p\underline{\bm{\sigma}}\left(\mathbf{x}_p\right)\cdot (\hat{\mathbf{D}}_{i,e}^v)^{-1} \underline{\bm{\sigma}}\left(\mathbf{x}_p\right) J_p w_p}{2E_{i,e}\sum_p J_p w_p}
=:
\dfrac{\h_{i,e}^{\bm \sigma}}{E_{i,e}}
\end{equation}
and

\begin{equation}
\bar{\Psi}_{ij,e}^{\bm \sigma,d}
=
\dfrac{\sum_p\underline{\bm{\tau}}\left(\mathbf{x}_p\right)\cdot (\hat{\mathbf{D}}_{ij,e}^d)^{-1} \underline{\bm{\tau}}\left(\mathbf{x}_p\right) J_p w_p}{2G_{ij,e}\sum_p J_p w_p}
=:
\dfrac{\h_{ij,e}^{\bm \sigma}}{G_{ij,e}}.
\end{equation}

The complementary \textit{ad-hoc} variables are

\begin{equation}
\h^{\bm \sigma}_{i,e} := \dfrac{\sum_p\underline{\bm{\sigma}}\left(\mathbf{x}_p\right)\cdot (\hat{\mathbf{D}}_{i,e}^v)^{-1} \underline{\bm{\sigma}}\left(\mathbf{x}_p\right) J_p w_p}{2\sum_p J_p w_p},
\end{equation}
and

\begin{equation}
\h^{\bm \sigma}_{ij,e} := \dfrac{\sum_p\underline{\bm{\tau}}\left(\mathbf{x}_p\right)\cdot (\hat{\mathbf{D}}_{ij,e}^d)^{-1} \underline{\bm{\tau}}\left(\mathbf{x}_p\right) J_p w_p}{2\sum_p J_p w_p}.
\end{equation}

On the one hand, the volumetric reduced matrices $(\hat{\mathbf{D}}_{i,e}^v)^{-1}$ are now direct to obtain, namely

\begin{equation}
[\hat{\mathbf{D}}^v_{1,k}] = \begin{bmatrix}
1 & 0 & 0 \\ -\nu_{12} & 0 & 0 \\ -\nu_{13} & 0 & 0
\end{bmatrix}, \;
[\hat{\mathbf{D}}^v_{2,k}] = \begin{bmatrix}
0 & -\nu_{21} & 0 \\ 0 & 1 & 0 \\ 0 & -\nu_{23} & 0
\end{bmatrix}, \;
[\hat{\mathbf{D}}^v_{3,k}] = \begin{bmatrix}
0 & 0 & -\nu_{31} \\ 0 & 0 & -\nu_{32} \\ 0 & 0 & 1
\end{bmatrix}.
\label{Eq:2_inverse_D_i^v}
\end{equation}

Therefore, the condition that isotropic materials meet it is likewise satisfied in the complementary formulation, that is, $\partial \h^{\bm \sigma}_{i,e} / \partial E_{i,e} = 0$ and $\partial \h^{\bm \sigma}_{ij,e} / \partial G_{ij,e} = 0$. This will imply a slightly simpler formulation that avoids computing derivatives of constitutive matrix terms.

On the other hand, for the deviatoric part, it turns out that the reduced matrices are the same as the direct formulation i.e.

\begin{equation}
(\hat{\mathbf{D}}_{ij,e}^d)^{-1} = \hat{\mathbf{D}}_{ij,e}^d,
\end{equation}
so the expressions for these matrices are the same as those stated in Equation~(\ref{Eq:2_D_ij_e_reduced}). Taking all into account, the total (complementary) energy is thus computed as

\begin{equation}
\W_e^{\bm \sigma}
=
\displaystyle\sum_{i=1}^3 \dfrac{\h^{\bm \sigma}_{i,e}}{E_{i,e}}\Omega_e
+ 
\displaystyle\sum_{i=1}^3 \sum_{j>i}^3 \dfrac{\h^{\bm \sigma}_{ij,e}}{G_{ij,e}} \Omega_e.
\label{Eq:2_W_e_c_decoupled_and_H}
\end{equation}

\subsection{Update formula derivation}

We now perform an extension starting from the isotropic algorithm formulation for direct and complementary energies. From Ben-Yelun et al. \cite{ben2023topology}, the minimization of the standard deviation $s_{\heps}(\bm \heps)$ of the isotropic strain-level variable $\bm \heps$ subjected to the equilibrium equation is performed applying a gradient-based scheme setting a step parameter $\eta_e := \heps_e N / k$, where $k$ is a (user-selected) update parameter to control the smoothness of the step. This yields the update formula

\begin{equation}
^{t+1}E_e = {}^{t}E_e\left(1 + \dfrac{^t\heps_e - {}^t\overline{\bm \heps}}{ s_{\bm \heps} k }\right),
\end{equation}
where $\overline{\bm \heps}$ is the average value of this mechanical property over all the elements in the domain.

We can extend this isotropic update parameter to the strain-based orthotropic formulation i.e., $\eta_{i,e}:=\heps_{i,e}N/k_i$ and $\eta_{ij,e} := \heps_{ij,e}N/k_{ij}$, ultimately arriving to

\begin{equation}
\begin{array}{c}
{^{t+1}}E_{i,e} =  \tensor[^t]{E}{_{i,e}} \left(1 + \tensor[^t]{\alpha}{_{i,e}} \right),\\[1ex]
{^{t+1}}G_{ij,e} = \tensor[^t]{G}{_{ij,e}} \left(1 + \tensor[^t]{\alpha}{_{ij,e}} \right),
\end{array}
\label{Eq:2_update_direct}
\end{equation}
with
\begin{equation}
\begin{array}{c}
\alpha_{i,e} = \dfrac{\heps_{i,e} - \overline{\bm \heps_i}}{s_{\bm \heps_i}k_i},\\
[3ex]
\alpha_{ij,e} = \dfrac{\heps_{ij,e} - \overline{\bm \heps_{ij}}}{ s_{\bm \heps_{ij}} k_{ij} }.
\end{array}
\label{Eq:2_alpha_definition}
\end{equation}

Therefore, update formulae to homogenize all six $\bm \heps_i$ and $\bm \heps_{ij}$ variables are derived computing the update parameters with the strains field information. Note that $k_i, k_{ij}$, are the (user-prescribed) modulating parameters and can be defined independently for every elastic property that is being optimized. For ease of notation  we write
\begin{equation}
    \bm{k} = \left\{\begin{array}{cccccc}
        k_1 & k_2 & k_3 & k_{12} & k_{13} & k_{23}
    \end{array} \right\}.
\end{equation}

Similarly, according to \cite{ben2023topology} and making use of previous extension, the update formulae for the stress-based orthotropic formulation are

\begin{equation}
\begin{array}{c}
^{t+1}E_{i,e} = \dfrac{^tE_{i,e}}{1 - \tensor[^t]{\alpha}{_{i,e}}}, \\[3ex]
^{t+1}G_{ij,e} = \dfrac{^tG_{ij,e}}{1 - \tensor[^t]{\alpha}{_{ij,e}}},
\end{array}
\label{Eq:2_update_complementary}
\end{equation}
where the update parameters $\alpha_{i, k}$ and $\alpha_{ij,e}$ are computed by means of Equation~\eqref{Eq:2_alpha_definition}, i.e., the strain-based case, since the strain homogenization is pursued. The difference lies in the way these variables are computed, making use of the following handy relations

\begin{equation}
\heps_{i,e} = \dfrac{\hsig_{i,e}}{E_{i,e}^2},
\quad
\heps_{ij,e} = \dfrac{\hsig_{ij,e}}{G_{ij,e}^2}.
\end{equation}

\subsection{Energy restrictions}
There are some constraints that have to be taken into account in order to prevent the optimization algorithm from creating a thermodynamically inconsistent material \cite{jones2018mechanics}. For isotropic materials, these are non-negative Young modulus $E>0$ and a bounded Poisson coefficient i.e. $-1 < \nu < 0.5$---these limits ensure the shear and bulk moduli to be positive, respectively.

The extrapolation to the orthotropic case of these restrictions are taken from Jones \cite{jones2018mechanics}, and are listed below:

\begin{equation}
\begin{array}{c}
\begin{array}{rl}
E_i \geq 0, & i = 1, 2, 3\\[2ex]
G_{ij} \geq 0, &  i=1,2,3; \; j>i\\[2ex]
\left|\nu_{ij} \right| < \sqrt{\dfrac{E_i}{E_j}}, & i=1,2,3; \; j>i
\end{array}\\[2ex]
\nu_{12} \nu_{23}\nu_{31} < \dfrac{1}{2}\left(1 - \nu_{12}\nu_{21} - \nu_{13}\nu_{31} - \nu_{23}\nu_{32}\right)
\end{array}
\label{Eq:2_energy_restrictions}
\end{equation}

\subsection{Problem statement}
With everything detailed above, the optimization problem might be stated as

\begin{equation}
\left\lbrace
\begin{array}{rl}
\min & \lbrace s_{\bm \heps_i}, s_{\bm \heps_{ij}}\rbrace, \quad i=1,2,3;\; j>i \\[2ex]
\textrm{s.t.} & \bm{K}(E_i, G_{ij}, \nu_{ij})\; \bm u = \bm f \\[2ex]
 & \begin{array}{ll}
 	0 \leq E_i \leq E_{\max}, & i = 1,2,3\\[2ex]
 	0 \leq G_{ij} \leq G_{\max}, &  i=1,2,3; \; j>i\\[2ex] 
	\left|\nu_{ij} \right| < \sqrt{\dfrac{E_i}{E_j}}, & i=1,2,3; \; j>i 
 \end{array}\\[2ex]
& \nu_{12} \nu_{23}\nu_{31} < \dfrac{1}{2}\left(1 - \nu_{12}\nu_{21} - \nu_{13}\nu_{31} - \nu_{23}\nu_{32}\right) 
\end{array}
\right.
\label{Eq:2_optimization_problem}
\end{equation}

All six standard deviations are minimized at once, by performing steps on the elastic properties likewise at once. This is achieved by separating every iteration into two steps: first to solve the equilibrium via FEM using the elastic properties from the previous iterations -- or the initial values, in case of iteration 0 -- thus obtaining the displacement, stress and strain distributions. Then the elastic properties are updated via Eqs.~(\ref{Eq:2_update_direct}) for the direct approach or (\ref{Eq:2_update_complementary}) if the complementary approach is used. This update is performed leaving the rest of the variables fixed.

The optimization is run until some convergence criterion is satisfied, setting a tolerance $\epsilon$ to this effect. Namely, the update is applied until the ratio between two successive computations of the elastic energy is lower than this tolerance i.e.

\begin{equation}
\dfrac{\left|^{t+1}\W^{\square} - {}^t\W^{\square}\right|}{^t\W^{\square}} \leq \epsilon.
\label{Eq:2_convergence}
\end{equation}

\section{Methodology}
\label{Sec:methodology}

In this section the specifics regarding the methods carried out in this paper are introduced. Therefore, the following subsections address the initial conditions that will be fed to the optimization algorithm, as well as the processes within the algorithm itself.  

%
%

\subsection{Boundary Conditions} 

Once the initial volume has been properly defined and meshed, the program will assign to each FE a set of initial properties, those being the six elastic properties considered in this work. The next step is now to prescribe BCs at the boundaries $\Gamma_D$ and $\Gamma_N$. Imposing continuous BCs along the cube is impossible due to the nature of FEs. Instead, the BCs must be prescribed on the nodes of the mesh. In order to select a set of nodes, several functions have been implemented with the purpose of selecting regions of the mesh and determining which nodes are contained within said regions. Having selected the nodes, it is possible to impose either displacements or forces on the nodes. Displacements would be measured in $\SI{}{\milli \meter}$, while force would be measured in $\SI{}{\newton}$. It must be noted that since tetrahedrons are solid FEs, it is not possible to impose rotations, since solid FEs do not have rotational DOFs.

\subsection{Optimization Algorithm}

With the domain and boundaries properly defined, it is now possible to commence the optimization. The stiffness update process is composed of three consecutive sub-processes. First of all, a base update will take place following Equations~\eqref{Eq:2_update_direct} or \eqref{Eq:2_update_complementary} for direct or complementary approaches, respectively. Secondly, some elements that have reached certain values are removed from the optimization. Finally, the resulting stiffness distribution is adjusted in order to fit the energy constraints of orthotropic materials.

\subsubsection{Base Update}

The base update process is based on the Equation (\ref{Eq:2_update_direct}) for the direct approach and Equation (\ref{Eq:2_update_complementary}) for the complementary approach. Additionally, we have proposed that the update variables $\tensor[^{t}]{\alpha}{_{i,e}}$ and $\tensor[^{t}]{\alpha}{_{ij,e}}$ shall take the form stated in Equation~(\ref{Eq:2_alpha_definition}). Recall that the parameters in $\bm{k}$ must be selected manually by the user and play a fundamental role on the results that will be attained, since they determine whether the update process is more or less aggressive. Several tests with different values of $\bm{k}$ have been performed in order to find an effective set of values for each scenario.

Furthermore, having three of the six Poisson ratios defined, lead to different forms of coming up with the three remaining. Therefore, two different ways of dealing with this computation have been implemented. The first method consists of fixing the values of $\nu_{12}, \nu_{23}$ and $\nu_{13}$, then assigning the corresponding values to $\nu_{21}, \nu_{32}$ and $\nu_{31}$ through the symmetry relation $\nu_{ij}/E_i = \nu_{ji} / E_j$.  Alternatively, we can impose a relationship such as $\nu_{ij} + \nu_{ik} = \kappa$ and then solve a system of equations to obtain $\nu_{ij}$ and $\nu_{ik}$.

Finally, once the FEs have been updated following the previous procedure, elastic properties that have reached values beyond their specified range will be updated to fit inside said range. The entire base update process is described in Algorithm \ref{alg:baseupd}.

\begin{algorithm}[p]
\caption{Base stiffness update.}
\label{alg:baseupd}
\begin{algorithmic}[1]
\Require ${\rm approach} = {\rm direct} \vee {\rm complementary}$
\For{$e = 1, \dots, {\rm number\_of\_FEs}$}
    \For{$i = 1, 2, 3$}
        \State $\alpha_{i,e} \gets \frac{\heps_{i,e}-\overline{\bm\heps_i}}{s_{\bm \heps_i}k_i}$
        \If{${\rm approach} = {\rm direct}$}
            \State $E_{i,e} \gets E_{i,e} \left(1 + \alpha_{i,e} \right)$
        \Else[${\rm approach} = {\rm complementary}$]
            \State $E_{i,e} \gets E_{i,e} \left(\frac{1}{1 - \alpha_{i,e}} \right)$
        \EndIf
        \If{$E_{i,e} > E_{\max}$}
            \State $E_{i,e} \gets E_{\max}$
        \ElsIf{$E_{i,e} < E_{\min}$}
            \State $E_{i,e} \gets E_{\min}$
        \EndIf
    \EndFor
    \For{$ij \in \{12, 23, 13\}$}
        \State $\alpha_{ij,e} \gets \frac{\heps_{ij,e}-\overline{\bm\heps_{ij}}}{s_{\bm \heps_{ij}}k_{ij}}$
        \If{${\rm approach} = {\rm direct}$}
            \State $G_{ij,e} \gets G_{ij,e} \left(1 + \alpha_{ij,e} \right)$
        \Else[${\rm approach} = {\rm complementary}$]
            \State $G_{ij,e} \gets G_{ij,e} \left(\frac{1}{1 - \alpha_{ij,e}} \right)$
        \EndIf
        \If{$G_{ij,e} > G_{\max}$}
            \State $G_{ij,e} \gets G_{\max}$
        \ElsIf{$G_{ij,e} < G_{\min}$}
            \State $G_{ij,e} \gets G_{\min}$
        \EndIf
    \EndFor
\EndFor
\end{algorithmic}
\end{algorithm}

\subsubsection{Element Removal}

In order to improve the performance of the algorithm and facilitate reaching convergence, certain FEs, or more specifically the elastic properties assigned to them, are removed from the optimization process and their value will be fixed when given lower and upper thresholds are surpassed. As to achieve those results, elastic properties that have reached values above the maximum allowed or below the minimum will not be updated in the following iterations.

\subsubsection{Fulfillment of Energy Constraints}

The values obtained after the base update may not be compatible with a real orthotropic material. Hence, they must be adjusted in order to satisfy the constraints introduced in Equation~(\ref{Eq:2_energy_restrictions}). The code will check in every FE whether the energy constraints are satisfied. If they are not, the values of the corresponding elastic properties must be altered. Two different alteration algorithms are used depending on the elastic energy associated to a FE compared to the average energy of all FEs. The stiffness of FEs with higher elastic energies are, in a first approximation, more relevant to the behaviour of the structure. Therefore, FEs that surpass a certain energy threshold will undergo a process that exclusively increases the values for the Young and shear moduli whenever possible, while the rest may be subject to a process that reduces some of the elastic properties in order to meet the required constraints.

Note that the constraints in Equations (\ref{Eq:2_energy_restrictions}.1) and (\ref{Eq:2_energy_restrictions}.2) do not need to be enforced since the base update prevents the Young and shear moduli from taking values below a certain minimum. Regarding the constraint in Equation (\ref{Eq:2_energy_restrictions}.3), it is always met when the second method of computing all the Poisson ratios is selected and as long as $\left\lvert \kappa \right\rvert < 2$.

This section has described in detail the procedures used to develop the object of this paper. This provides the understanding required to introduced the obtained results and the conclusions that can be inferred from them in the following sections.

\section{Numerical examples}
\label{Sec:results}

This section presents the results that have been obtained using the procedures explained in Section~\ref{Sec:methodology}. To this end, four different load cases are studied. For each of them, figures with the achieved stiffness distribution are presented, as well as a figure indicating which FEs have the greatest elastic energy and the numerical value of the overall compliance of the structure, which is the main factor at evaluating the effectiveness of each method. The results are computed using the direct approach for the isotropic case, and both the direct and complementary approaches for the orthotropic case. It has been deemed unnecessary to show the results for the isotropic case using both approaches since the results obtained are very similar between them, specially when large values of $\bm{k}$ are selected---see Ben-Yelun et al. for more details \cite{ben2023topology}.

\begin{figure}[htb]
    \centering
    \includegraphics[width=1.0\textwidth]{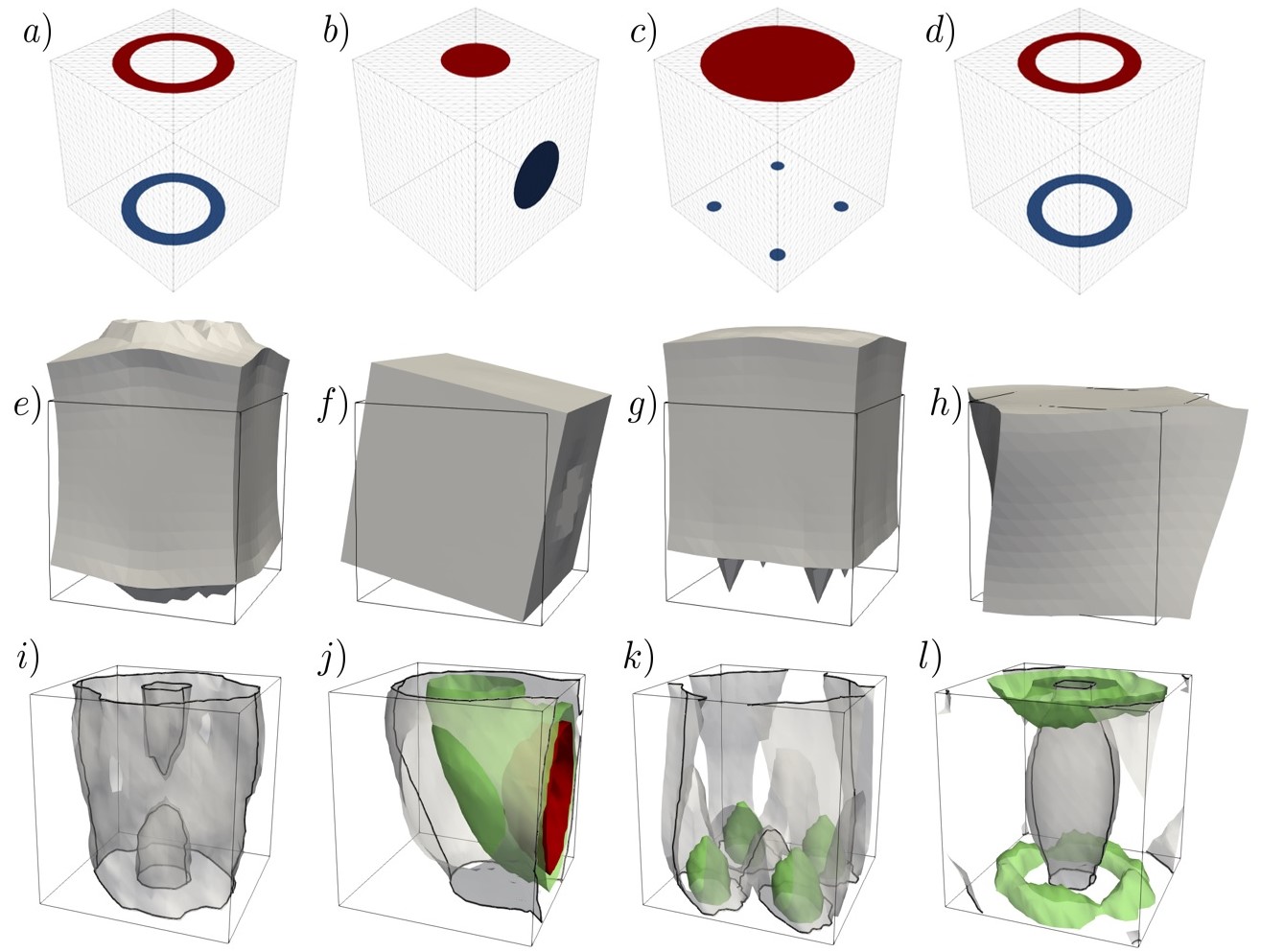}
    \caption{The boundary conditions of the four cases studied in this section, denoted as tube (a), elbow (b), chair (c) and torsion (d). The first row are the boundary conditions (a-d), the second row is the initial deformed configuration of the cube with different magnifications (e-h), and the third row is the isosurface of elastic energy for the optimized sample using the orthotropic calculation with the complementary approach (i-l). In grey the isosurface of  0.8 $\mu$J/mm$^3$, in green the one of 8 $\mu$J/mm$^3$ and in red the 80 $\mu$J/mm$^3$ isosurface.}
    \label{fig:num1}
\end{figure}

Figure \ref{fig:num1} shows the four cases studied. All use the same initial cube and properties, but with different boundary conditions. The top row shows the boundary conditions, the middle row the deformation imposed in each case, with different magnifications, and the bottom row shows the isosurfaces with a constant elastic energy density for the orthotropic complementary calculation of each case.

These results will be contrasted between each other and against the results from a traditional SIMP optimization from the commercial program \textsc{Altair} \textsc{OptiStruct}\textsuperscript{\textregistered}.

\subsection{Hyperparameters}

For all the load cases, the generator cube has a side length of $\SI{100}{\milli\meter}$ and there are 12 subdivisions per side, for a total of 8640 FEs. The initial and limit values assigned assigned to the FEs for the isotropic case are:

\begin{gather*}
   \tensor[^0]{E}{} = \SI{200}{\giga\pascal}; \quad \nu = 0.33, \\
   E_{\max} = \SI{500}{\giga\pascal}; \quad E_{\min} = \SI{50}{\giga\pascal},
\end{gather*}
and, for the orthotropic case,
\begin{gather*}
\tensor[^0]{E}{_{i}} = \SI{200}{\giga\pascal}; \quad \tensor[^0]{G}{_{ij}} = \SI{75.19}{\giga\pascal}; \quad \nu_{ij} = 0.33, \\
E_{i,\max} = \SI{500}{\giga\pascal}; \quad E_{i,\min} = \SI{50}{\giga\pascal}, \\
G_{ij,\max} = \SI{75.19}{\giga\pascal}; \quad G_{ij,\min} = \SI{7.519}{\giga\pascal}; \quad i,j = 1,2,3; \enspace j > i.
\end{gather*}

In all cases the defined values for $\nu_{ij}$ remain constant during the optimization process (i.e. $\nu_{12}$, $\nu_{13}$ and $\nu_{23}$, recall that the others are updated along with the elastic properties within the optimization loop). The parameters in $\bm{k}$ will be defined individually for each load case. Finally, the tolerance established for the converging condition in Equation~(\ref{Eq:2_convergence}) takes the value $\epsilon = 0.01$.

In the case of \textsc{OptiStruct}\textsuperscript{\textregistered} the generator cube has the same dimensions but it is meshed in a different fashion, with elements of average size \SI{5}{\milli\meter}. The SIMP method deals with a material with a maximum Young modulus $E_{\max} = \SI{500}{\giga\pascal}$ using $p=1$. The penalization parameter is set to 1 since this solution might stiffen some elements using less volume fraction of them, thus achieving a lower compliance in comparison with analyses with $p>1$. Furthermore, using a finer mesh allows the problem to reach an slightly lower compliance. Both of these adjustments benefit the final compliance obtained with the SIMP method, thus enabling a more conservative focus for the methods presented in this paper. Finally, the volume fraction is set to $\hat{V} = V/V_{\max} = 0.5$ as a trade-off solution, since lower volume fractions do not let the domain to be completely (topologically) connected, and higher fractions tend to over-stiffen the final structure, thus achieving a pointless low-compliance solution. Also, the final shape obtained using this volume fraction is similar to the ones obtained using the other methods.

\subsection{Load Case 1: Tube}

The first load case, named Tube, is characterized by normal stresses along the $z-$axis without many stresses on any other direction. For this scenario, all the parameters in $\bm{k}$ are set to 15. The imposed boundary conditions to simulate this load case are displayed in \figurename~\ref{fig:num1}a. This load case is defined by two circular crowns on the top and bottom faces of the generator cube. The bottom crown (blue) represents the Dirichlet boundary $\Gamma_D$, where all displacements are set to $\SI{0}{\milli\meter}$. The top crown (red) represents part of the Neumann boundary $\Gamma_N$, where a force of $\SI{10000}{\newton}$ in the direction of the $z-$axis has been distributed among all the contained nodes. The remaining surface of the cube also belongs to the boundary $\Gamma_N$, since the prescribed force on all the nodes contained within is equal to $\SI{0}{\newton}$.


\subsubsection{Isotropic Case}

\begin{figure}[htb]
    \centering
    \includegraphics[width=0.6\textwidth]{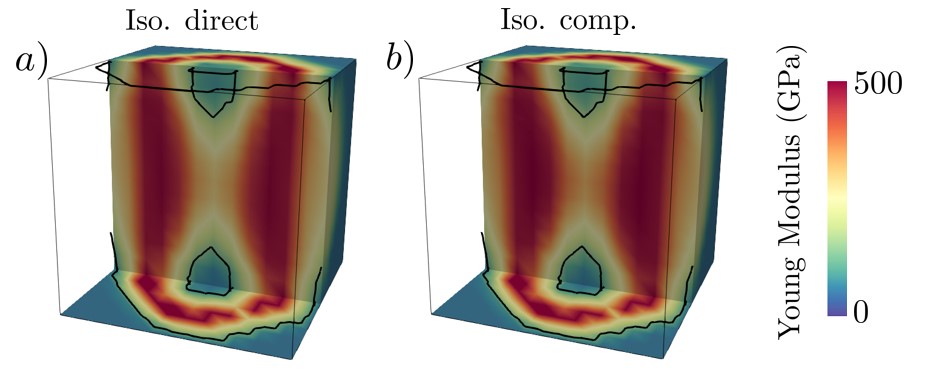}
    \caption{Optimal Young moduli distribution for the direct and complementary isotropic calculations for the tube case.}
    \label{fig:isotropic1}
\end{figure}

The isotropic case resulted in a compliance of $\SI{2.23}{\newton\!\,\milli\meter}$ after 21 iterations and provided the stiffness distributions in Figure \ref{fig:isotropic1}.

The obtained results seem very reasonable since both the Young modulus and the elastic energy are the greatest in a tube-like volume between the two circular crowns, as one might have expected in the first place. It must be noted that only FEs with elastic energies greater than a certain value relative to the maximum elastic energy within the generator cube have been represented. The quasi-transparent FEs have surpassed a less exigent threshold, while the opaque ones have surpassed a more exigent threshold.

This case will serve as the benchmark to evaluate the two following results, which are related to the orthotropic case.

\subsubsection{Orthotropic Case}

For this load case, and the successive orthotropic analyses as well, only three out of all the six elastic properties are commented i.e. the most relevant ones. In the Tube load case, the elastic properties of interest are $E_{x}$, $E_{z}$ and $G_{yz}$.

\paragraph{Direct Approach}

\begin{figure}[htb]
    \centering
    \includegraphics[width=0.8\textwidth]{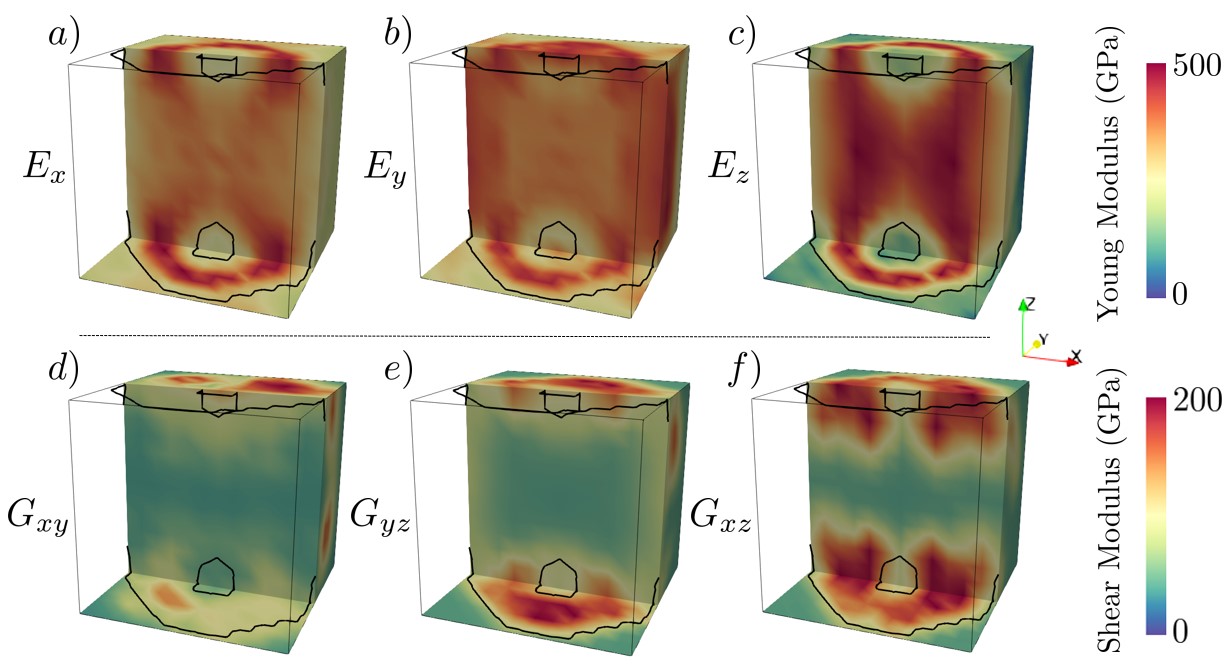}
    \caption{Optimal distribution for the tube considering the direct approach in the orthotropic case.}
    \label{fig:direct1}
\end{figure}

The orthotropic case using the direct approach resulted in a compliance of $\SI{2.03}{\newton\!\,\milli\meter}$ after 14 iterations and provided the stiffness distributions in Figure \ref{fig:direct1}.

Similarly to the isotropic case, $E_{z}$ develops into a tube-like structure between the circular crowns. On the other hand, $E_{x}$ -- analogous to $E_{y}$ -- is greater next to the circular crowns and on the faces perpendicular to the $y-$axis, most likely to minimize the compliance associated to the Poisson effect. Lastly, $G_{yz}$ -- analogous to $G_{xz}$ -- is greatest next to both of the circular crowns, presenting a significant directionality.

The FEs with the most elastic energy are organized in a very similar structure to that of the isotropic case, although more concentrated next to the faces. This rearrangement may explain the better local optimum reached in this orthotropic analysis.

\paragraph{Complementary Approach}

\begin{figure}[htb]
    \centering
    \includegraphics[width=0.8\textwidth]{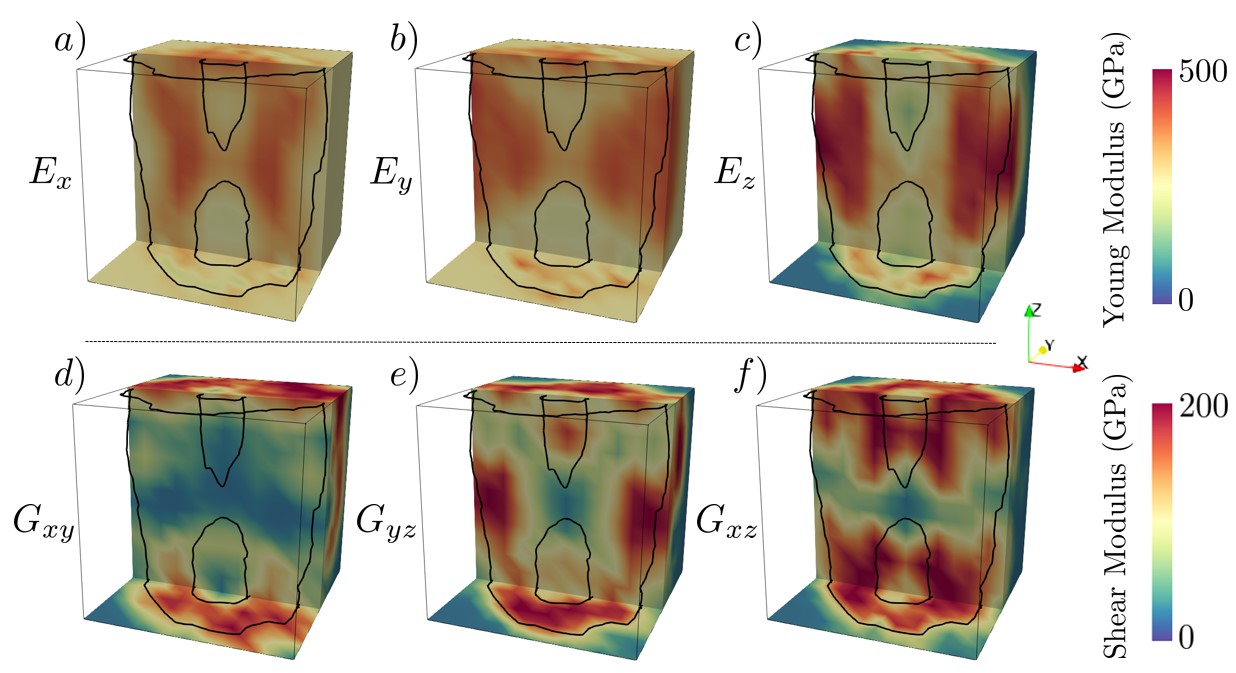}
    \caption{Optimal distribution for the tube considering the complementary approach in the orthotropic case.}
    \label{fig:complementary1}
\end{figure}

The orthotropic case using the complementary approach resulted in a compliance of $\SI{0.897}{\newton\!\,\milli\meter}$ after 28 iterations and provided the stiffness distributions in Figure \ref{fig:complementary1}. 

In this case $E_{z}$ does not take the shape of a full tube, losing stiffness next to the top and bottom faces and reaching greater values in the rest of them. This can most likely be attributed to the fact that those areas have been stiffened by $G_{yz}$ instead. Additionally, it seems like $E_{x}$ just ignores the tube-like expected shape in order to deal with the Poisson effect as efficiently as possible.

Once again, the FEs with the most elastic energy are organized in a very similar structure to that of the isotropic case, although wider.

\subsubsection{Comparison}

It seems like the isotropic and orthotropic direct cases reach a similar solution judging by the similar stiffness distribution regarding $E_{z}$ and the compliance attained being just slightly smaller for the isotropic case.

However, the orthotropic complementary case reaches a much lower compliance with a less expected \textit{a priori} distribution. Lacking experimental results performed on real materials, it is impossible to determine whether such a result would be reproducible on an actual structure, but it holds significant value as a result regardless.

\begin{figure}[htb]
    \centering
    \includegraphics[width=0.7\textwidth]{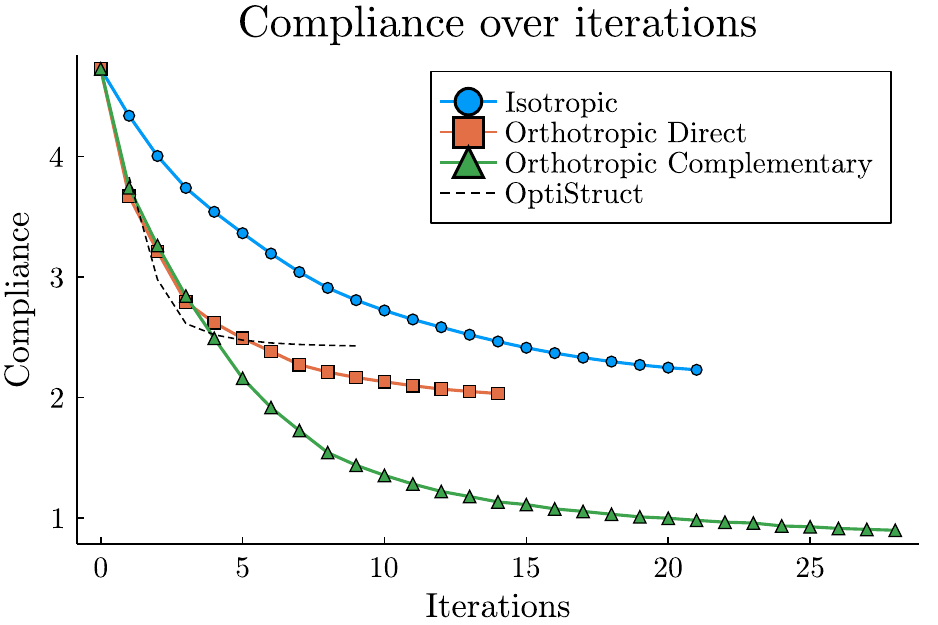}
    \caption{Evolution of the compliance over the number of iterations for Tube.}
    \label{fig:evolution1}
\end{figure}

The evolution of the compliance with the number of iterations performed for the previous methods, as well as \textsc{OptiStruct}\textsuperscript{\textregistered}, can be found in Figure \ref{fig:evolution1}. These results show that every approach used outperforms the traditional optimization using a SIMP method with the (hyper-)parameters highlighted before.

\begin{figure}[htb]
    \centering
    \includegraphics[width=1.0\textwidth]{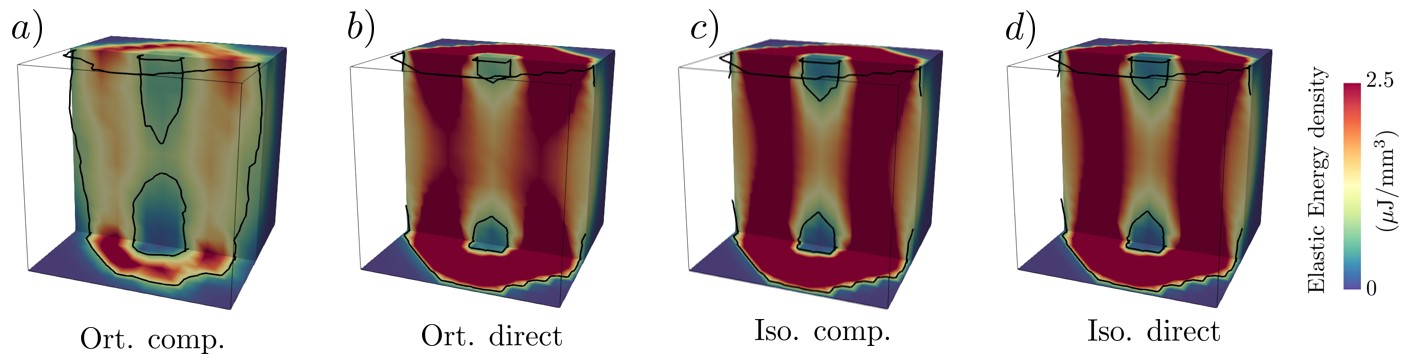}
    \caption{Elastic energy density distribution in the optimal configuration of the four calculation cases; orthotropic complementary (a), orthotropic direct (b), isotropic complementary (c) and isotropic direct (d).}
    \label{fig:enercomptube}
\end{figure}

Figure \ref{fig:enercomptube} shows the comparison of the final elastic energy density distribution considering the 4 configurations used in the calculation (i.e. isotropic direct, isotropic complementary, orthotropic direct, orthotropic complementary). Both the isotropic approaches reach practically the same result, but the orthotropic cases can release the elastic energy density from the middle region of the tube. It means that the orthotropic methods are more efficient, and even more the complementary approach.

\subsection{Load Case 2: Elbow}

The second load case, named Elbow, is characterized by normal and shear stresses along the $x$ and $z-$axis. For this scenario, all the parameters in $\bm{k}$ are set to 20. The representation of this load case is depicted in \figurename~\ref{fig:num1}b. This load case is defined by two circles on the face perpendicular to the $x-$axis and the top face: the former (blue) represents the Dirichlet boundary $\Gamma_D$, where all displacements are set to $\SI{0}{\milli\meter}$. The latter (red) represents part of the Neumann boundary $\Gamma_N$, where a force of $\SI{10000}{\newton}$ in the direction of the $z-$axis has been distributed among all the contained nodes. Analogous to the first load case, the remaining surface of the cube also belongs to the boundary $\Gamma_N$, since the prescribed force on all the nodes contained within is equal to $\SI{0}{\newton}$.

\subsubsection{Isotropic Case}

\begin{figure}[htb]
    \centering
    \includegraphics[width=0.6\textwidth]{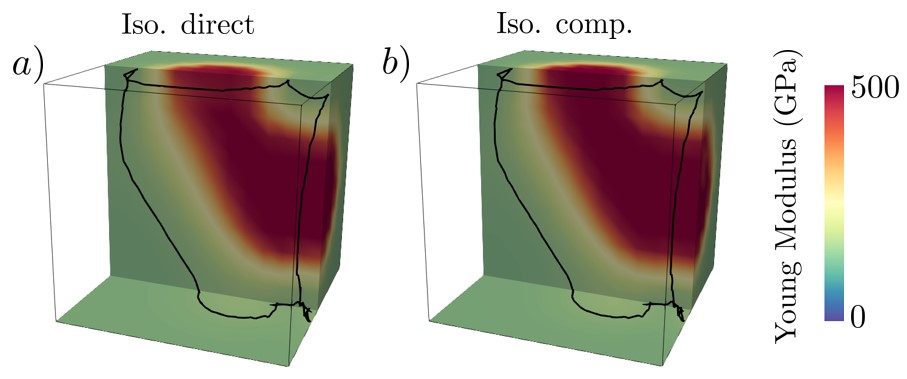}
    \caption{Optimal Young moduli distribution for the direct and complementary isotropic calculations for the elbow case.}
    \label{fig:isotropic2}
\end{figure}

The isotropic case resulted in a compliance of $\SI{23.2}{\newton\!\,\milli\meter}$ after 15 iterations and provided the stiffness distributions in Figure \ref{fig:isotropic2}.

The Young modulus distribution seems to reproduce a structure similar to an elbow joint, joining both the boundary conditions. The elastic energy distribution adopts a similar shape, where the FEs with the greatest energy are grouped next to the Dirichlet boundary. This case is compared with the two following (orthotropic) analyses.

\subsubsection{Orthotropic Case}

The relevant elastic properties, which are therefore highlighted, are $E_{x}$, $E_{z}$ and $G_{xz}$ since these are the directions involved in this load case.

\paragraph{Direct Approach}

\begin{figure}[htb]
    \centering
    \includegraphics[width=0.8\textwidth]{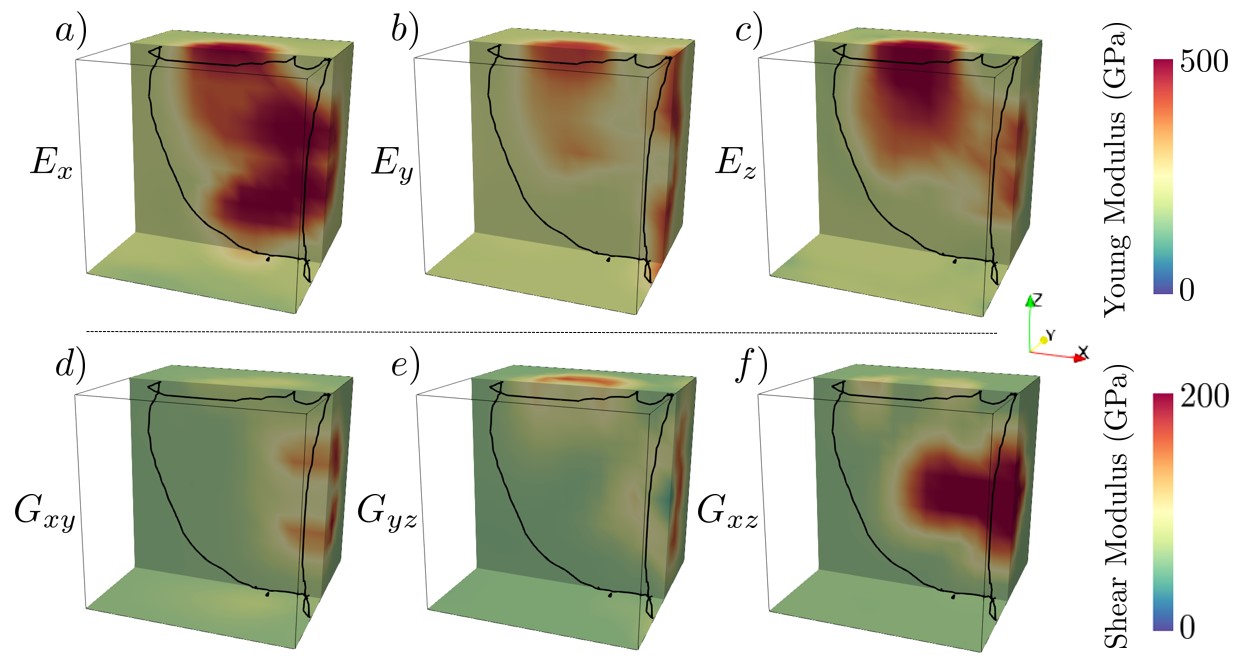}
    \caption{Optimal distribution for the elbow considering the direct approach in the orthotropic case.}
    \label{fig:direct2}
\end{figure}

The orthotropic case using the direct approach resulted in a compliance of $\SI{22.4}{\newton\!\,\milli\meter}$ after 12 iterations and provided the stiffness distributions in Figure \ref{fig:direct2}. 

The Young modulus along $x-$axis $E_{x}$ develops in a structure similar to an elbow joint -- i.e. similarly to the isotropic case -- although more prominent near the Dirichlet boundary. On the other hand, $E_{z}$ seems to only reproduce the part of the elbow closer to the Neumann boundary, while $E_{z}$ has the same behaviour next to the Dirichlet boundary. The FEs with the most elastic energy are organized in a near identical structure to that of the isotropic case.

\paragraph{Complementary Approach}

\begin{figure}[htb]
    \centering
    \includegraphics[width=0.8\textwidth]{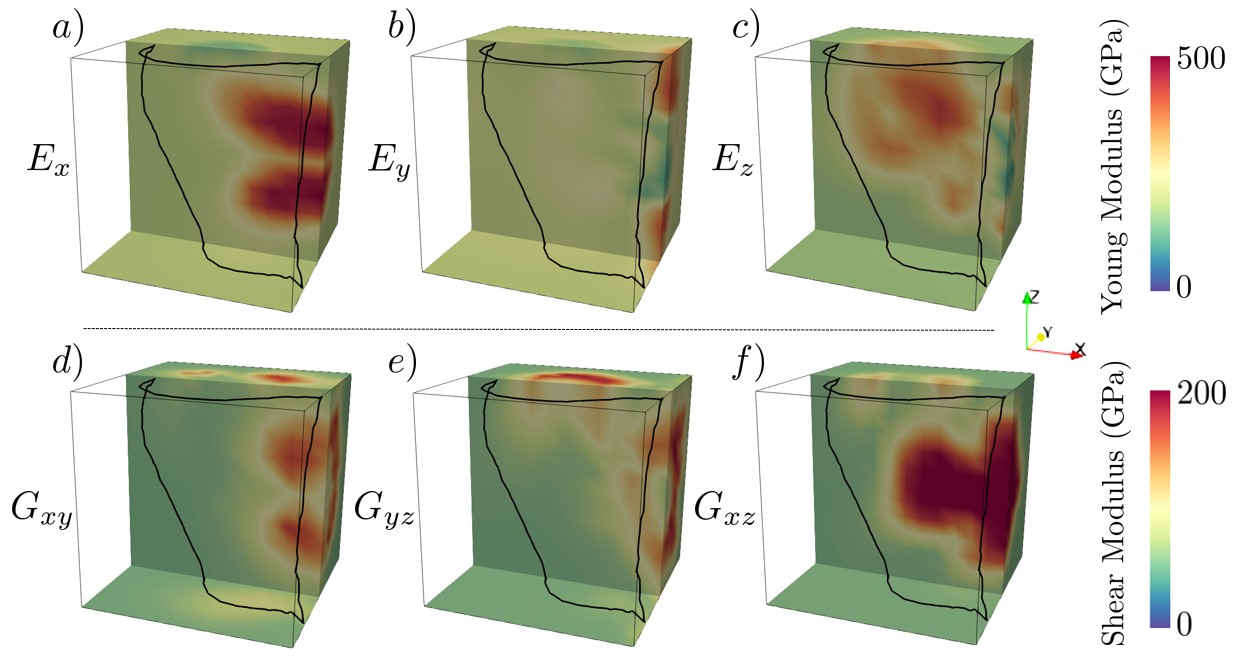}
    \caption{Optimal distribution for the elbow considering the complementary approach in the orthotropic case.}
    \label{fig:complementary2}
\end{figure}

The orthotropic case using the complementary approach resulted in a compliance of $\SI{22.8}{\newton\!\,\milli\meter}$ after 14 iterations and provided the stiffness distributions in Figure \ref{fig:complementary2}. 

In this case, the stiffness and energy distributions reached are almost identical to those of the direct approach. The main difference being $E_{x}$ is now only significant next to the Dirichlet boundary, more precisely the top and bottom parts, indicating that it is most likely trying to stiffen the areas where there is a large bending moment.

\subsubsection{Comparison}

It seems like all cases reach similar solutions judging by the stiffness and energy distributions and the final compliance being very similar. The main point of interest of this load case is how the orthotropic case is capable of discerning which elastic property is the adequate one to stiffen a certain area and minimize the compliance of the structure as a result. 

\begin{figure}[htb]
    \centering
    \includegraphics[width=0.7\textwidth]{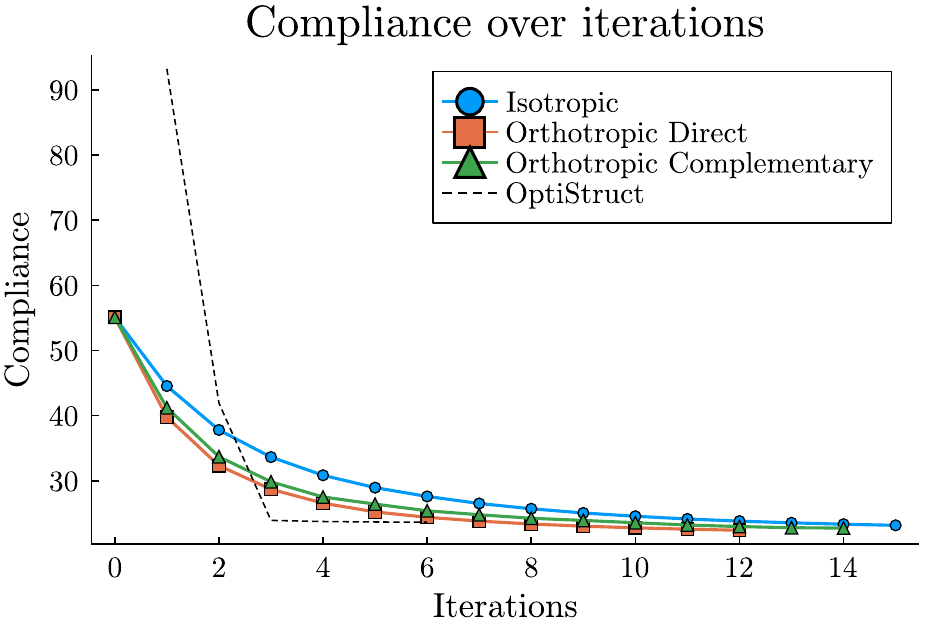}
    \caption{Evolution of the compliance over the number of iterations for Elbow.}
    \label{fig:evolution2}
\end{figure}

The evolution of the compliance with the number of iterations performed for the previous methods, as well as \textsc{OptiStruct}\textsuperscript{\textregistered}, can be found in Figure \ref{fig:evolution2}. These results show that every approach used matches the performance of  the traditional optimization using a SIMP method.

\begin{figure}[htb]
    \centering
    \includegraphics[width=1.0\textwidth]{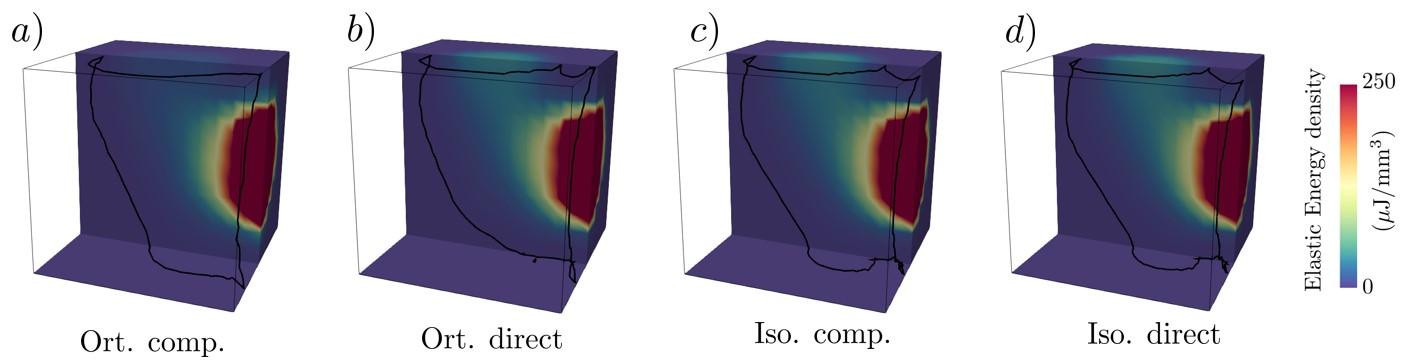}
    \caption{Elastic energy density distribution in the optimal configuration of the four calculation cases; orthotropic complementary (a), orthotropic direct (b), isotropic complementary (c) and isotropic direct (d).}
    \label{fig:enercompelbow}
\end{figure}

Figure \ref{fig:enercompelbow} shows that the elastic energy density distribution in the sample is practically de same in all the calculation cases, which agrees with the behaviour of Figure \ref{fig:evolution2}.

As it might be observed in \figurename~\ref{fig:evolution2}, our model has succeeded to equalize the SIMP method, yet no improvement is observed. One of the reasons that may explain this result is the coexistence of loads of two different nature i.e. FEs that are subjected to axial and shear loads. Therefore, despite the fact that the orthotropic case enables the optimization of both longitudinal and shear stiffness separately, it is difficult to meet a trade-off solution since the effects that they mutually produce might not be decoupled, and the algorithm finds it difficult to perform big steps along the optimization---note the smooth changes in compliance in the three curves.

\subsection{Load Case 3: Chair}

The third load case, named Chair, is characterized by normal stresses along the $z-$axis and the presence of significant shear stresses. For this scenario, all the parameters in $\bm{k}$ are set to 30. As seen in Figure \ref{fig:num1}c, this load case is defined by four small circles on the bottom face and a bigger circle on the top face. The small circles (blue) represent the Dirichlet boundary $\Gamma_D$, where all displacements are set to $\SI{0}{\milli\meter}$. The bigger circle (red) represents part of the Neumann boundary $\Gamma_N$, where a force of $\SI{10000}{\newton}$ in the direction of the $z-$axis has been distributed among all the contained nodes.

\subsubsection{Isotropic Case}

\begin{figure}[htb]
    \centering
    \includegraphics[width=0.6\textwidth]{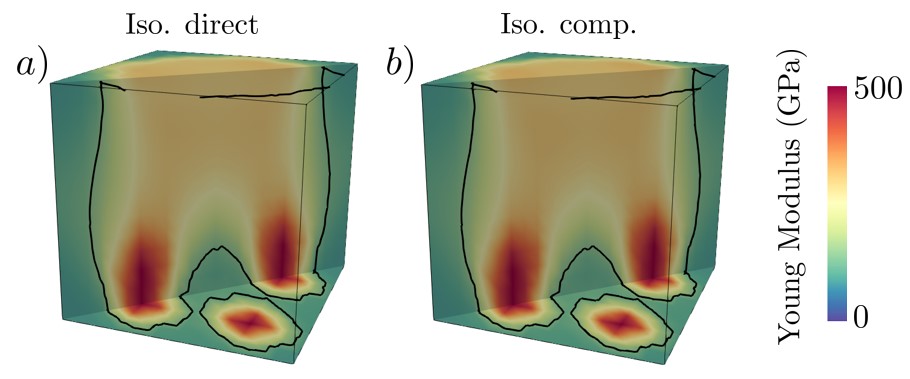}
    \caption{Optimal Young moduli distribution for the direct and complementary isotropic calculations for the chair case.}
    \label{fig:isotropic3}
\end{figure}

The isotropic case resulted in a compliance of $\SI{4.40}{\newton\!\,\milli\meter}$ after 20 iterations and provided the stiffness and elastic energy distributions in Figure \ref{fig:isotropic3}.

The Young modulus distribution seems to reproduce a structure reminiscent of a stool, being the most stiff near the surfaces conforming the Dirichlet boundary. The elastic energy distribution is nearly identical, both in shape and magnitude.

\subsubsection{Orthotropic Case}

In this case, the elastic properties of interest are $E_{z}$, $G_{xy}$ and $G_{yz}$, conclusions on the latter are valid for $G_{xz}$ as well.

\paragraph{Direct Approach}

\begin{figure}[htb]
    \centering
    \includegraphics[width=0.8\textwidth]{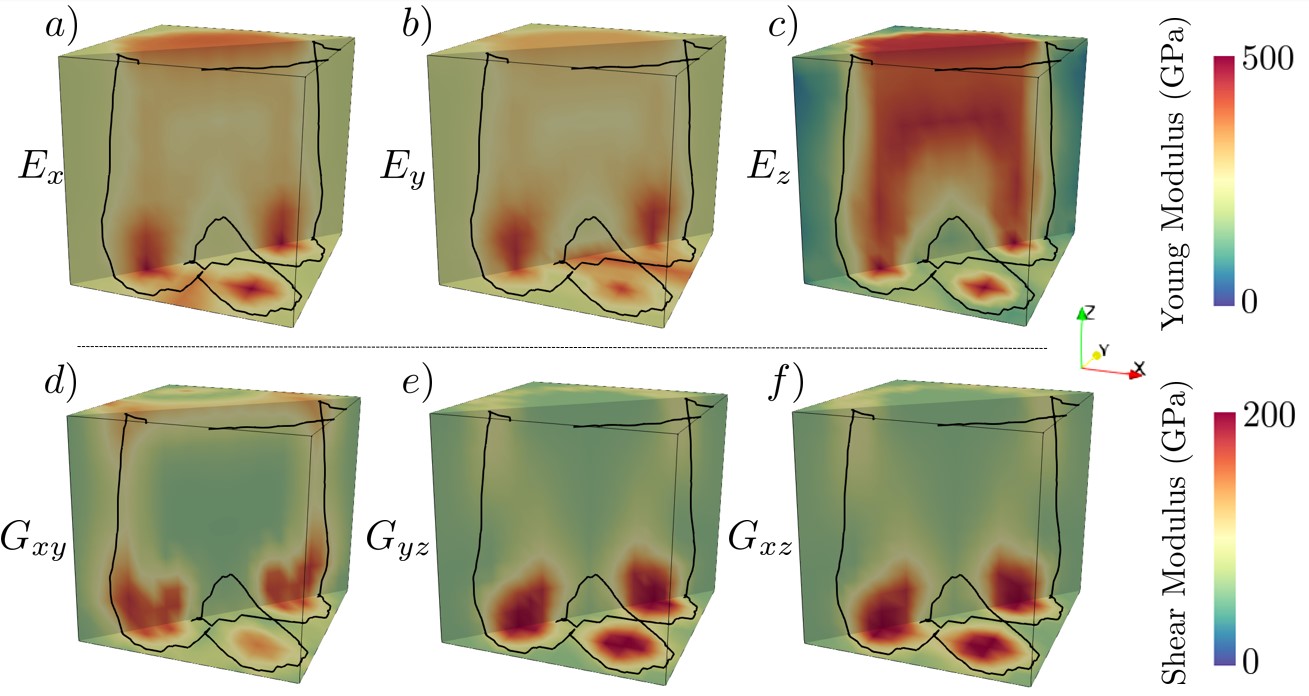}
    \caption{Optimal distribution for the chair considering the direct approach in the orthotropic case.}
    \label{fig:direct3}
\end{figure}

The orthotropic case using the direct approach resulted in a compliance of $\SI{3.75}{\newton\!\,\milli\meter}$ after 16 iterations and provided the stiffness and elastic energy distributions in Figure \ref{fig:direct3}. 

The result in $E_{z}$ is analogous to the isotropic computation since it seems to reproduce a structure reminiscent of a stool, although in this case the stiffness is distributed more evenly throughout the structure. However, the `legs' of the stool are still the most stiff parts since both $G_{xy}$ and $G_{yz}$ are the most prominent here. Regarding the elastic energy distribution along the FEs, those with higher values are organized in a near identical structure to that of the isotropic case.

\paragraph{Complementary Approach}

\begin{figure}[htb]
    \centering
    \includegraphics[width=0.8\textwidth]{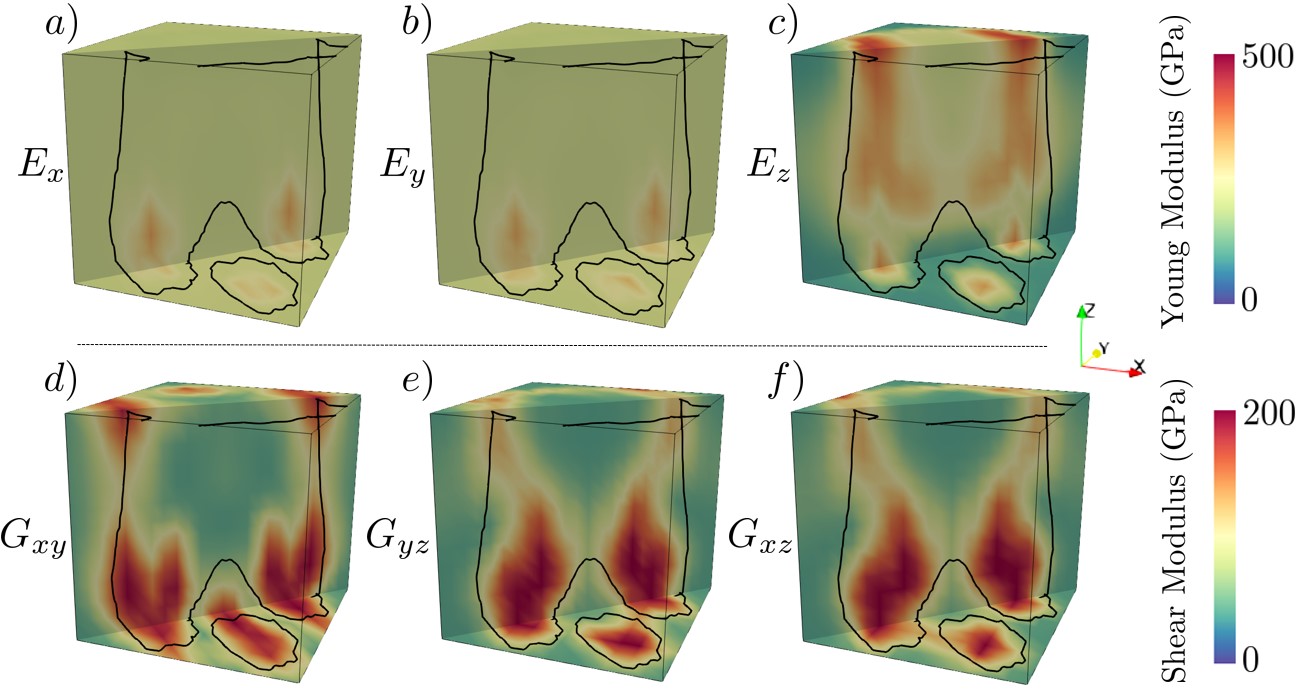}
    \caption{Optimal distribution for the chair considering the complementary approach in the orthotropic case.}
    \label{fig:complementary3}
\end{figure}

The orthotropic case using the complementary approach resulted in a compliance of $\SI{1.87}{\newton\!\,\milli\meter}$ after 30 iterations and provided the stiffness and elastic energy distributions in Figure \ref{fig:complementary3}. 

In this case the structure developed by $E_{z}$ is still reminiscent of a stool, but it is significantly less stiff than in the previous two cases. Similarly to the direct approach, both $G_{xy}$ and $G_{yz}$ reach their higher values near the Dirichlet boundary, but in this case $G_{yz}$ high values extend all the way up to the Neumann boundary following the contour of the stool. Additionally, it seems some areas have been over-stiffened, specially by both the shear moduli. Nevertheless, the FEs with the most elastic energy are still organized in a very similar fashion to the previous two cases.

\subsubsection{Comparison}

For this load case both the isotropic case and the orthotropic case reach qualitatively similar results for the stiffness distribution, although the latter manages to attain a significantly lower compliance.

\begin{figure}[htb]
    \centering
    \includegraphics[width=0.7\textwidth]{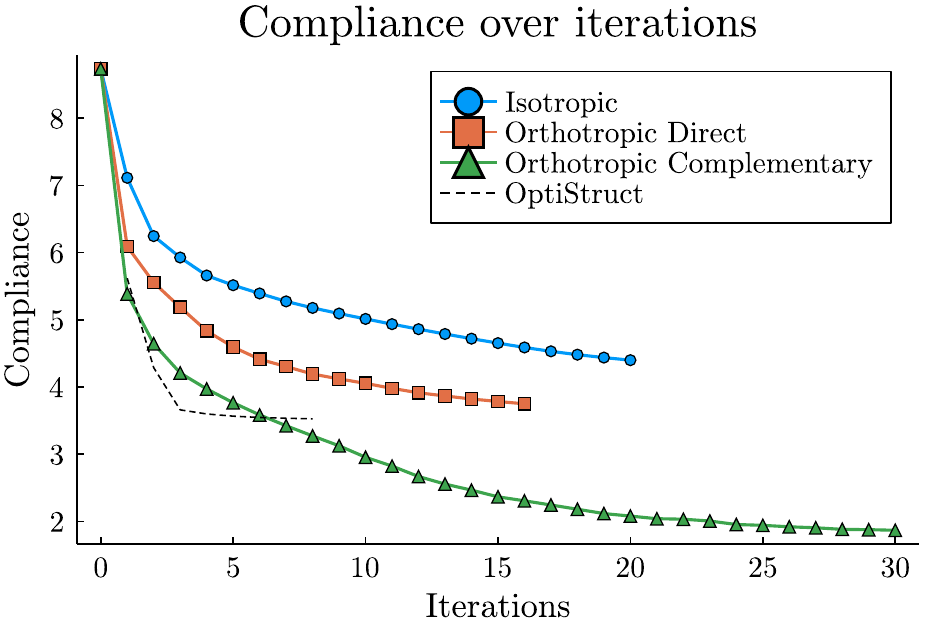}
    \caption{Evolution of the compliance over the number of iterations for Chair.}
    \label{fig:evolution3}
\end{figure}

However, just like in the first load case, the orthotropic complementary case succeeds in reaching a much lower compliance with a less intuitive distribution. As it has already been stated, this is an interesting result, deserving of further analysis in experimental testing. The evolution of the compliance with the number of iterations performed for the previous methods, as well as \textsc{OptiStruct}\textsuperscript{\textregistered}, can be found in Figure \ref{fig:evolution3}. 

\begin{figure}[htb]
    \centering
    \includegraphics[width=1.0\textwidth]{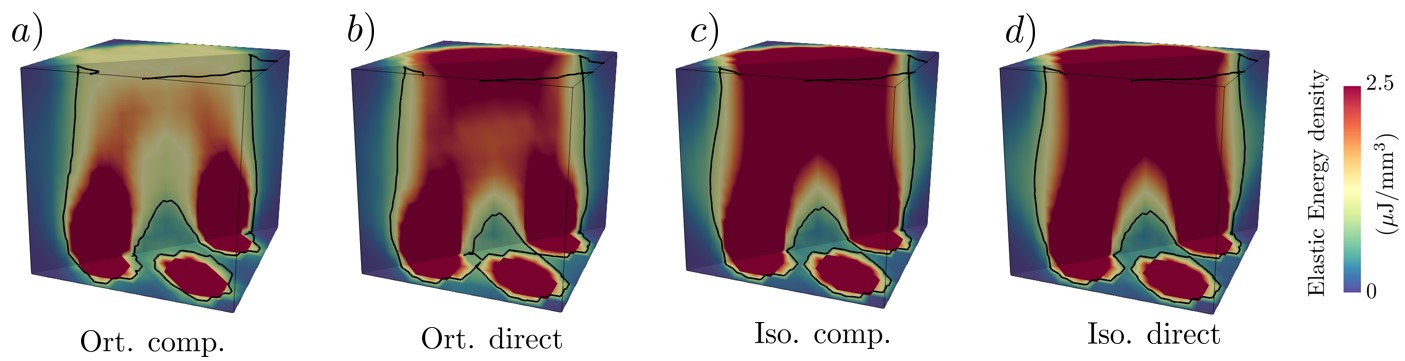}
    \caption{Elastic energy density distribution in the optimal configuration of the four calculation cases; orthotropic complementary (a), orthotropic direct (b), isotropic complementary (c) and isotropic direct (d).}
    \label{fig:enercompchair}
\end{figure}

Figure \ref{fig:evolution3} shows a clear improvement given by the orthotropic approach, even more with the complementary case. This is also clear in Figure \ref{fig:enercompchair}, where the isotropic cases have a high density of elastic energy in all the domain. And the orthotropic cases are capable to release this elastic energy density in the top part of the sample, with a big difference between the orthotropic complementary approach and the others. 

The results displayed in \figurename~\ref{fig:evolution3} show significant disparity between the compliance obtained by the different methods, where the orthotropic case using a complementary approach manages to outperform the rest (including SIMP) by a significant margin.

\subsection{Load Case 4: Torsion}

The fourth load case, named Torsion, is characterized by primarily shear stresses without any significant normal stresses in any direction. For this scenario, all the parameters in $\bm{k}$ are set to 20. Boundary conditions for this computation are represented in Figure~\ref{fig:num1}d. This load case is defined by two circular crowns on the top and bottom faces of the generator cube. The bottom crown (blue) represents the Dirichlet boundary $\Gamma_D$, where all displacements are set to $\SI{0}{\milli\meter}$. The top crown (red) represents part of the Neumann boundary $\Gamma_N$, where a force of $\SI{10000}{\newton}$ has been distributed among all the contained nodes. The direction of the force vector $\mathbf{f}_{\bm k}$ at each node is perpendicular to the vector that goes from the central point of the top face to that node $\mathbf{r}_{\bm{O,k}}$, and the sense is such that the cross product of the aforementioned vector times the force vector i.e. the torque is positive along the $z-$axis, $\left(\mathbf{r}_{\bm{O,k}} \times \mathbf{f}_{\bm k} \right) \cdot \mathbf{u}_{\bm z} > 0$.

\subsubsection{Isotropic Case}

\begin{figure}[htb]
    \centering
    \includegraphics[width=0.6\textwidth]{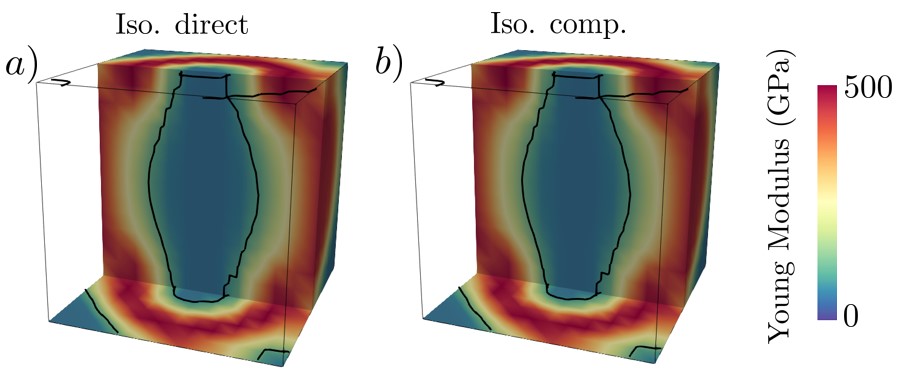}
    \caption{Optimal Young moduli distribution for the direct and complementary isotropic calculations for the torsion case.}
    \label{fig:isotropic4}
\end{figure}

The isotropic case resulted in a compliance of $\SI{3.81}{\newton\!\,\milli\meter}$ after 23 iterations and provided the stiffness and elastic energy distributions in Figure \ref{fig:isotropic4}.

The Young modulus acquires the greatest values near the two boundary conditions and in the outermost parts of the cube. Additionally, a clear directionality can be observed when comparing the faces perpendicular to the $x-$axis and those perpendicular to the $y-$axis. This phenomenon occurs because this is a pure shear stress case, so the principal directions for stresses and strains appear $\ang{45}$, thus the directionality in the stiffness distribution. Contrary to that, the innermost past has reached the minimum possible stiffness values. The elastic energy distribution adopts a similar shape, where the FEs with the greatest energy are grouped next to both of the boundaries.

\subsubsection{Orthotropic Case}

In this case, the elastic properties of interest are $E_{x}$ -- conclusions of which apply to $E_{y}$ as well --  $G_{xz}$ and $G_{yz}$. The distributions $E_{z}$ and $G_{xy}$ are of lesser significance due to their low contribution to the elastic energy of the structure.

\paragraph{Direct Approach}

\begin{figure}[htb]
    \centering
    \includegraphics[width=0.8\textwidth]{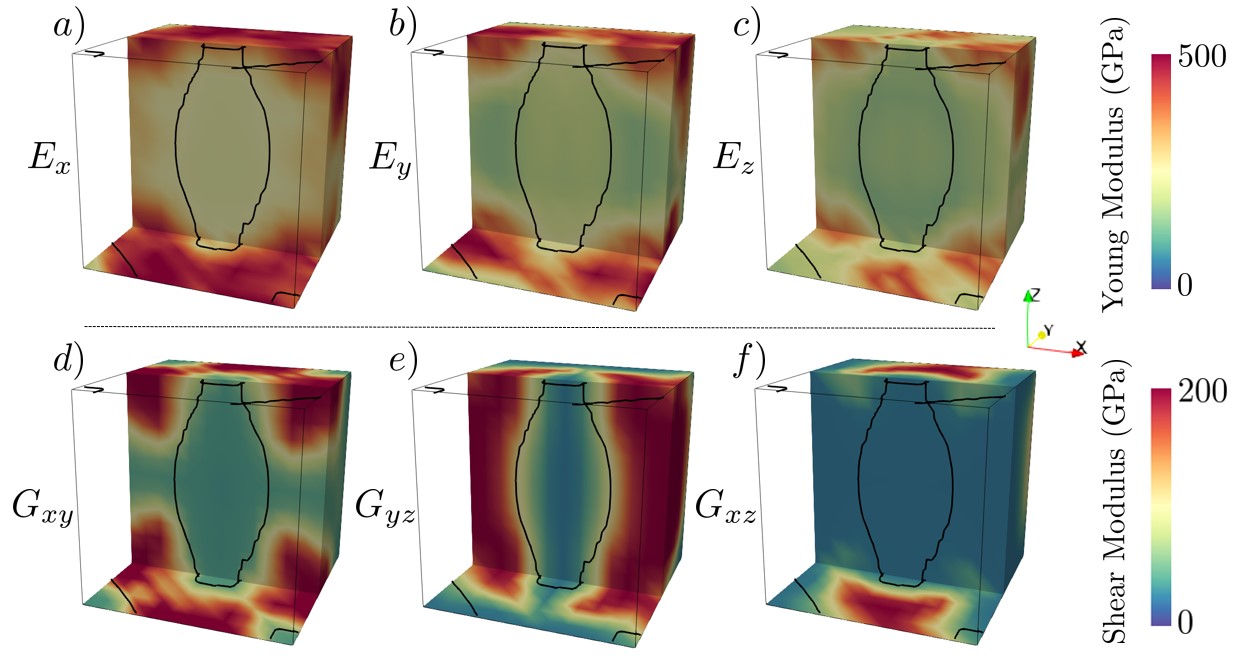}
    \caption{Optimal distribution for the torsion considering the direct approach in the orthotropic case.}
    \label{fig:direct4}
\end{figure}

The orthotropic case using the direct approach resulted in a compliance of $\SI{3.11}{\newton\!\,\milli\meter}$ after 23 iterations and provided the stiffness and elastic energy distributions in Figure \ref{fig:direct4}. 

In this case it seems that the shear moduli $G_{xz}$ and $G_{yz}$ replicate the same structure as in the isotropic case, splitting the total area between the two, each modulus being responsible of stiffening two opposite faces of the cube. However, it is also apparent that the Young modulus $E_{x}$ has severely over-stiffened the cube while developing into a strange structure.

The FEs with the most elastic energy are organized in a very similar structure to that of the isotropic case, although more concentrated next to the top and bottom faces.

\paragraph{Complementary Approach}

\begin{figure}[htb]
    \centering
    \includegraphics[width=0.8\textwidth]{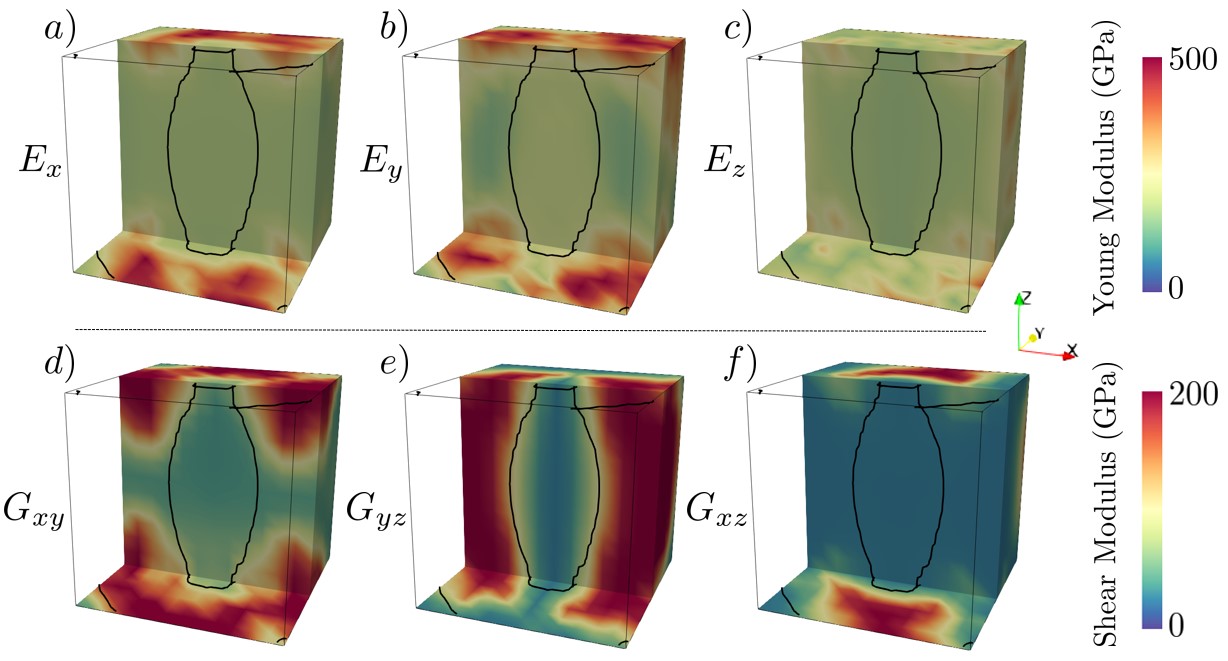}
    \caption{Optimal distribution for the torsion considering the complementary approach in the orthotropic case.}
    \label{fig:complementary4}
\end{figure}

The orthotropic case using the complementary approach resulted in a compliance of $\SI{2.87}{\newton\!\,\milli\meter}$ after 23 iterations and provided the stiffness and elastic energy distributions in Figure \ref{fig:complementary4}. 

For this case a very similar stiffness distribution to that of the direct approach has been reached by $G_{xz}$ and $G_{yz}$. Furthermore, $E_{x}$ does not present the same strange behaviour and over-stiffening as in the previous case. It is instead localised near the top and bottom faces with a very clear directionality. This leads to a better optimum solution reached.

Once again, the FEs with the most elastic energy are organized in a very similar structure to that of the isotropic case, although slightly more concentrated next to the top and bottom faces.

\subsubsection{Comparison}

This load case clearly represents the advantages of orthotropic materials, since the stiffness distribution can be clearly divided in two different regions depending on which shear modulus is required in that area and the compliance reached is significantly and consistently lower. Furthermore, it suggest that the complementary approach may be more suitable to avoid over-stiffening for the less relevant elastic properties.

\begin{figure}[htb]
    \centering
    \includegraphics[width=0.7\textwidth]{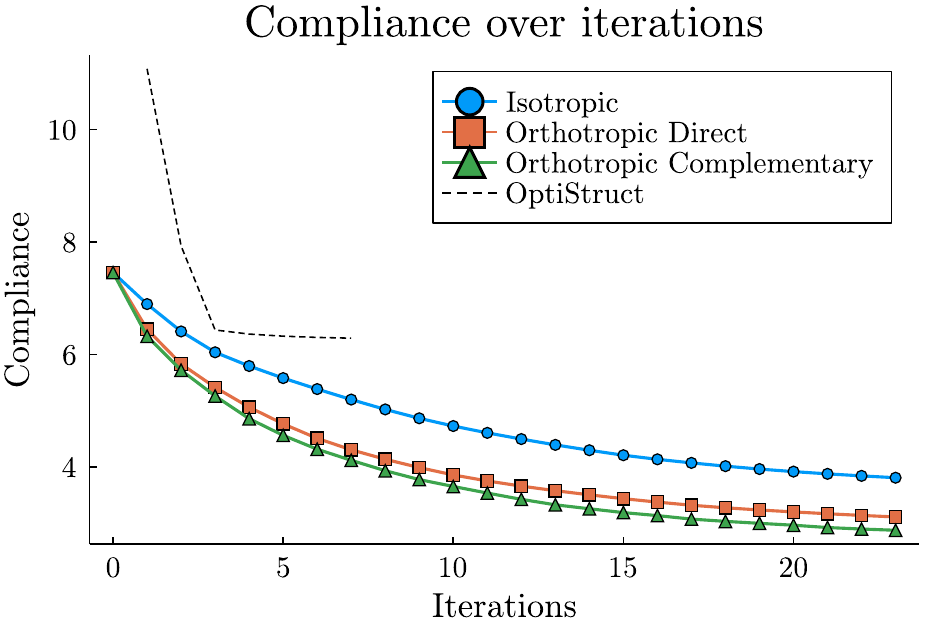}
    \caption{Evolution of the compliance over the number of iterations for Torsion.}
    \label{fig:evolution4}
\end{figure}

The evolution of the compliance with the number of iterations performed for the previous methods, as well as \textsc{OptiStruct}\textsuperscript{\textregistered}, can be found in Figure \ref{fig:evolution4}. As one can observe all the methods presented in this paper, and specially the orthotropic case with a complementary approach, show a clear advantage when compared to SIMP.

\begin{figure}[htb]
    \centering
    \includegraphics[width=1.0\textwidth]{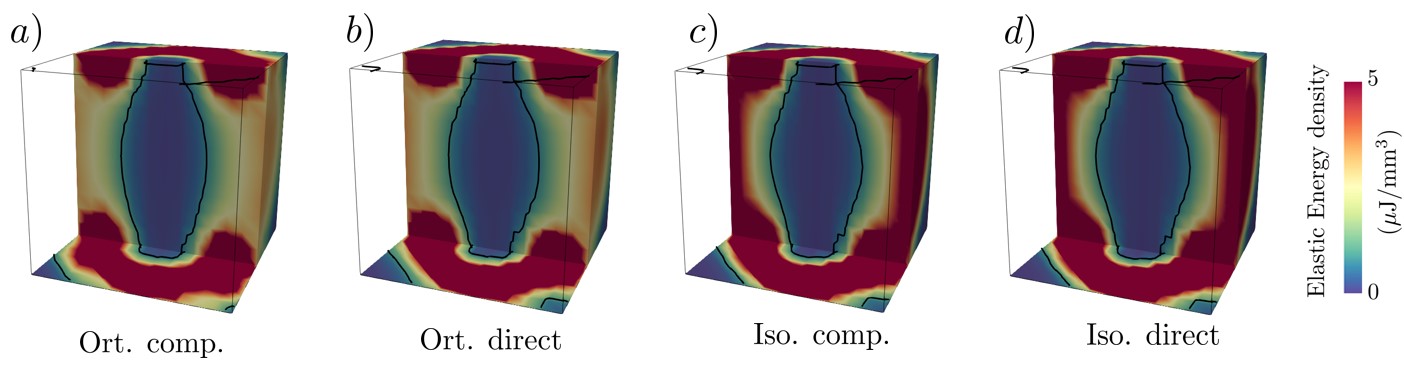}
    \caption{Elastic energy density distribution in the optimal configuration of the four calculation cases; orthotropic complementary (a), orthotropic direct (b), isotropic complementary (c) and isotropic direct (d).}
    \label{fig:enertorsion}
\end{figure}

Figure \ref{fig:enertorsion} shows a lower elastic energy density in the orthotropic cases, but very similar between them. As well, the elastic energy density distribution is similar for the isotropic cases. It agrees with the overall behaviour shown in Figure \ref{fig:evolution4}.

\subsection{Comparison of all load cases}

Table~\ref{tab:compliance} displays the overall results for the compliance obtained in all cases.

\begin{table}[htb]
    \centering
    \begin{tabular}{ I l I *{4}{C{1.8cm} I} }
        \bottomrule[\heavyrulewidth]
        \rowcolor{Gray}\textbf{Compliance [N\,mm]} & Tube & Elbow & Chair & Torsion  \\
        \bottomrule[\heavyrulewidth]
        \textsc{OptiStruct} & $2.43$ & $23.7$ & $3.53$ & $6.29$ \\
        \hline
        Isotropic & $2.23$ & $23.2$ & $4.40$ & $3.81$ \\
        \hline
        Orthotropic Direct & $2.03$ & $\mathbf{22.4}$ & $3.75$ & $3.11$ \\
        \hline
        Orthotropic Complementary & $\mathbf{0.897}$ & $22.8$ & $\mathbf{1.87}$ & $\mathbf{2.87}$ \\
        \bottomrule[\heavyrulewidth]
    \end{tabular}
    \caption{Final compliance obtained for every load case and optimization method.}
    \label{tab:compliance}
\end{table}

In nearly all cases, both the orthotropic results seem to be superior to those attained with the isotropic case or using \textsc{OptiStruct}\textsuperscript{\textregistered}. The only instance where this statement does not hold true is for the third load case (Elbow) when using the direct approach for the orthotropic case. Regardless of that, the orthotropic case either with a direct or a complementary approach manages to acquire the lowest compliance in all load cases.

Using these criteria and taking into account that the values obtained have been set up to be on the conservative side, it seems clear that using stiffness based methods may prove to be advantageous compared to a traditional SIMP in the scenarios herein presented.

\section{Concluding remarks}
\label{Sec:concluding_remarks}

We have successfully developed a formulation for linear orthotropic materials within a topology optimization framework through the development of an algorithm that homogenizes the strain level by performing an optimization on the elastic parameters of each FE, rather than using density-based methods which operate on an intermediate variable such as the volume of said FE. Only one parameter (or set of parameters with same value) $\bm{k}$ has to be fixed, and the energy constraints are ensured not to be violated. This paper represents an extension of a previous stiffness-based work valid for isotropic materials \cite{saucedo2023updated,ben2023topology}. The proposed generalization to orthotropic materials poses several advantages, most notably a higher versatility than its isotropic counterpart given that the design space is six times larger for the former---orthotropic materials present 9 elastic properties, however in this work we only consider 6 (longitudinal and shear moduli), fixing three Poisson ratios. Therefore, all of this leads to a procedure which is able to outperform the isotropic formulation. It also improves the results provided by an unpenalised SIMP implemented in the commercial software \textsc{OptiStruct}\textsuperscript{\textregistered}, with a similar usage of material. Note that the pursued objectives are different: whereas SIMP minimizes the compliance under a volume constraint, our method aims at minimizing the standard deviation of the strain variables $\bm\heps$. However, it is worth noting the similarities between the final shapes.

As has been highlighted, two approaches have been proposed to this end: the so-called direct, whose formulation is based on strains, and the complementary approach, developed within a stresses-based framework. In both methods, a first volumetric-deviatoric energy split is performed, then those terms are likewise split such that the contributions to energy and the effects of each of the six elastic properties are decoupled, thus enabling the update of all properties at the same time in every iteration. To this effect, 6 \textit{ad-hoc} update parameters $\alpha$ are defined through the strain elastic energy density variables of these decoupled terms, in a way that the convergence to local minimum of compliance is ensured, given by the application of a gradient-based scheme in the standard deviation of such variables. We have proved that the complementary formulation is more explicit than the direct -- as was originally highlighted by Suzuki and Kikuchi \cite{suzuki1991homogenization}, and Díaz and Bendsoe \cite{diaz1992shape} -- since complementary energy becomes more explicit in the elastic parameters, thus yielding more optimal results for the load cases herein presented. 

Regarding the improvements in results provided by this method, it has been observed for all the load cases that the orthotropic formulation compared to its analogous isotropic approach  achieves lower compliance structures with quite similar final shapes in terms of elastic energy. This is due to the previously commented versatility in the design space. This higher flexibility, result of a wider design space, enables new ways of distributing the energy contribution across the domain -- as Figures \ref{fig:enercomptube}, \ref{fig:enercompelbow}, \ref{fig:enercompchair} and \ref{fig:enertorsion} show -- hence the improvement regarding the isotropic case. Moreover, by improving this case results, we automatically outperform SIMP as well, since the isotropic formulation at least equalled this method having the same material usage \cite{ben2023topology}. Specially noteworthy is the fourth load case (torsion) in which a 54\% improvement of the compliance is achieved in comparable computation time using the complementary orthotropic optimization with respect to the SIMP method. We recall again that the pursued objectives are different.

New horizons with more structural meaningful objective functions that do not rely on constraint impositions can be exploited through proper structural optimization frameworks. This will bring in turn the ability to face the simulation, design, or optimization challenges that new families of materials are requiring, such as functionally graded materials or mechanical metamaterials. Our method aims at this objective.

\section{Acknowledgements}
\noindent\begin{minipage}[c]{.15\linewidth}
\includegraphics[width=\linewidth]{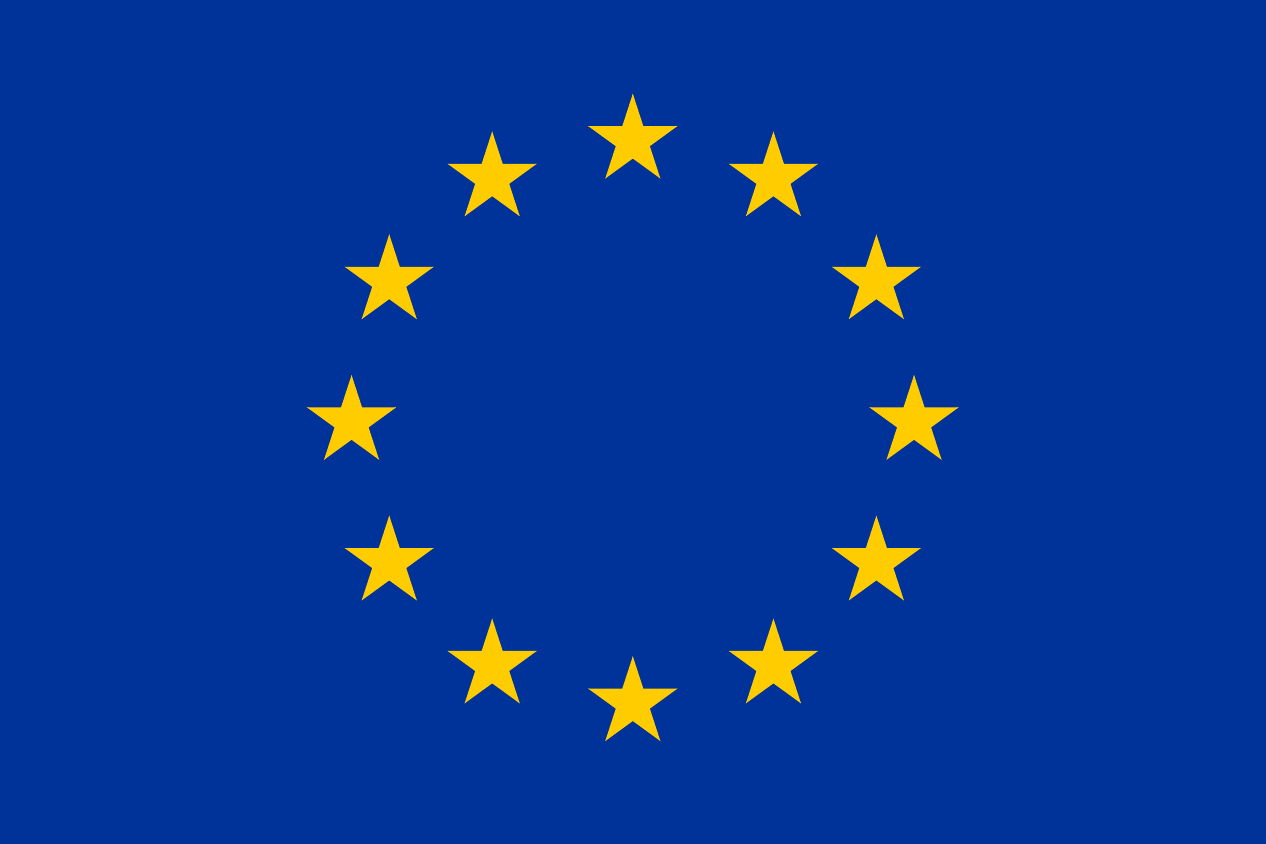}
\end{minipage}\hfill
\begin{minipage}[c]{.8\linewidth}
This project has received funding from the European Union's Horizon 2020 research and innovation programme under the Marie Skłodowska-Curie Grant Agreement No. 101007815.
\end{minipage}

\addcontentsline{toc}{section}{References}
\bibliographystyle{unsrt-custom}
\bibliography{references.bib}

\end{document}